\documentclass[12pt,reqno]{article}
\usepackage{amssymb}
\usepackage{amsmath}
\usepackage{amsthm}
\usepackage{amsfonts}
\usepackage{url}
\usepackage{hyperref}
\usepackage{breakurl}
\usepackage{rotating}
\usepackage{adjustbox}
\def\modd#1 #2{#1\ \mbox{\rm (mod}\ #2\mbox{\rm )}}
\usepackage{pdflscape}
\usepackage{microtype}
\newcommand{\comment}[1]{}
\usepackage{caption}
\newcommand{\BE}{\begin{equation}}
\newcommand{\EE}{\end{equation}}

\theoremstyle{plain}
\newtheorem{theorem}{Theorem}

\theoremstyle{definition}

\begin{document}
\title{Classical length-5 pattern-avoiding permutations}

\author{Nathan Clisby\\
Department of Mathematics,\\
Swinburne University of Technology,\\
Hawthorn, Vic. 3122, Australia\\
\href{}{\tt nclisby@swin.edu.au},\\
Andrew R Conway\\
Fairfield, Vic. 3078, Australia\\
\href{}{\tt andrewpermutations5@greatcactus.org},\\
Anthony J Guttmann\\
School of Mathematics and Statistics\\
The University of Melbourne\\
Vic. 3010, Australia\\
\href{mailto:guttmann@unimelb.edu.au}{\tt guttmann@unimelb.edu.au}
\and
Yuma Inoue\\Google Japan, SHIBUYA STREAM,\\
 3-21-3 Shibuya, Shibuya-ku,\\ Tokyo 150-0002, Japan\\
\href{mailto:guttmann@unimelb.edu.au}{\tt yumai@google.com}}

\date{}

\maketitle
\abstract
We have made a systematic numerical study of the 16 Wilf classes of length-5 classical pattern-avoiding permutations from their generating function coefficients. We have extended the number of known coefficients
in fourteen of the sixteen classes. Careful analysis, including sequence extension, has allowed us to estimate the growth constant of all classes, and in some cases to
estimate the sub-dominant power-law term associated with the exponential growth.

In six of the sixteen classes  we find the familiar power-law behaviour, so that the coefficients behave like $s_n \sim C \cdot  \mu^n \cdot n^g,$ while in the remaining
ten  cases we find a stretched exponential as the most likely sub-dominant term, so that the coefficients behave like $s_n \sim C \cdot \mu^n \cdot \mu_1^{n^\sigma} \cdot n^g,$ where $0 < \sigma < 1.$ We have also classified the 120 possible permutations into the 16 distinct classes.

We give compelling numerical evidence, and in one case a proof, that all 16 Wilf-class generating function coefficients can be represented as moments of a non-negative measure on $[0,\infty).$
Such sequences are known as {\em Stieltjes moment sequences}. They have a number of nice properties, such as log-convexity, which can be used to provide quite strong rigorous lower bounds.

Stronger bounds still can be established  under plausible monotonicity assumptions about the terms in the continued-fraction expansion of the generating functions implied by the Stieltjes property.
In this way we provide strong (non-rigorous) lower bounds to the growth constants, which are sometimes within a few percent of the exact value.

\section{Introduction}
\label{introduction}
Let $\pi$ be a permutation on $[n]$ and $\tau$ be a permutation on $[k].$ $\tau$ is said to occur as a {\em pattern} in $\pi$ if for some sub-sequence of $\pi$ of length $k$ all the elements of the sub-sequence occur in the same relative order as do the elements of $\pi.$ If the permutation $\tau$ does not occur in $\pi,$ then this is said to be a {\em pattern-avoiding permutation} or PAP.

Let $s_n(\tau)$ denote the number of permutations of $[n]$ that avoid the pattern $\tau.$ Stanley and Wilf conjectured, and Marcus and Tardos \cite{MT04} subsequently proved, that for any pattern $\tau$ in $[k]$ the limit $\lim_{n \to \infty} s_n(\tau)^{\frac{1}{n}}$ exists and is finite. This means that the number of PAPs grows exponentially with $n,$ whereas of course the number of permutations of $n$ grows factorially. 

There are 6 possible permutations of length three, and the number of permutations of length $n$ avoiding any of these 6 patterns is given precisely by $C_n = \frac{1}{n+1} \binom{2n}{n}\sim \frac{4^n}{\sqrt{\pi n^3}},$ where $C_n$ denotes the $n^{th}$ Catalan number. That is to say, all 6 possible patterns have the same exponential growth-rate as PAPs. Alternatively expressed, there is only one Wilf class for length-3 PAPs.

For length-4 PAPs there are three Wilf classes. Typical representatives of the three classes are $1234,$ $1342$ and $1324.$ The generating function for the first two classes is known. 

In the first case \cite{G90} the generating function is D-finite, satisfying a third-order linear, homogeneous ODE, and $s_n(1234) \sim \frac{81\sqrt{3}\cdot 9^n}{16\pi \cdot n^4}.$ 

In the second case \cite{B97} the generating function is algebraic, and $s_n(1342) \sim \frac{64\cdot 8^n}{243\sqrt{\pi} \cdot n^{5/2}}.$ 

The third case has not been solved, but extensive numerical work \cite{CGZ18} suggests that $s_n(1324) \sim C \cdot \mu^n \cdot \mu_1^{\sqrt{n}} \cdot n^g,$ where $\mu \approx 11.598$ (and possibly $9+3\sqrt{3}/2$ exactly), $\mu_1 \approx 0.040,$ and $g \approx -1.1.$ The appearance of the sub-dominant term $\mu_1^{\sqrt{n}}$ is referred to as a {\em stretched exponential} term.

Thus these three Wilf classes have generating functions that are D-finite, algebraic, and (almost certainly) non-D-finite respectively.

For the 16 Wilf classes of length-5 PAPs there is only one, $Av(12345),$ for which the generating function is known \cite{BEGM20}. It is D-finite. In one other case, $Av(31245),$ the growth constant is known \cite{B05}, and for $Av(53421)$ the growth constant can be expressed in terms of that of $Av(1324)$ \cite{B07}, which has been estimated to some accuracy in \cite{CGZ18}. These three known growth constant results can all be obtained from Theorem 4.2 in \cite{B07}, alternatively proved as Theorem 3.2 in \cite{APV19}.

We wrote a general purpose program to generate the coefficients, which is very efficient, but does have significant memory demands. It is memory, rather than time, that has limited the length of the series we can generate. There is considerable variation in both time and memory requirements for different Wilf classes. For example, for $Av(43251)$ we obtained the series to O$(x^{26})$ in 2 hours 50 mins of CPU time, using a single core, and 212 GB of memory. By contrast,  for $Av(52341)$ we obtained the series to O$(x^{23})$ in 15 hours of CPU time, using a single core, and 600 GB of memory.  To obtain the next coefficient would require more than 2TB of memory. In this way, we obtained series for all Wilf classes of various lengths from O$(x^{23})$ to O$(x^{27})$. 

This almost doubles the length of the available series in 14 of the 16 Wilf classes.  As mentioned, $Av(12345)$ is completely solvable, so an arbitrary number of coefficients is available, and for $Av(31245)$ a special purpose program has been written giving 38 coefficients \cite{B19}. Our program is described in the next section.

In Sec. 3 we give a brief description of Stieltjes moment sequences, describing the properties that we need in this work. In Sec. 4 we give a little more detail about what is known about classical patterns of length 5, and give a table classifying all 120 permutations into the 16 Wilf classes, based on a combination of known symmetries and direct enumeration. In Sec. 5 we describe the general principles of series analysis as needed here, in Sec. 6 we discuss the ratio method, and in Sec. 7 we discuss the extension of the ideas of the previous two sections to the analysis of stretched-exponential singularities. In Sec. 8 we describe the method of differential approximants, which is what is needed in Sec. 9, where we discuss the idea of series extension. More precisely, we are able to obtain many more {\em approximate} coefficients than we have exact coefficients to a sufficient degree of accuracy that we can apply the previously described ratio methods and their extensions.

Using all these techniques, we analyse the 16 Wilf classes in Sec. 10. Sec. 11 comprises a brief discussion and conclusion.

One {\em caveat} we would like to make is that the only singularity types we are considering
are pure power-law and stretched exponentials, motivated of course by the fact that these are the only singularity types we have encountered for PAPs of shorter length.
And indeed, it is clear that many of the series behave in one or other of these ways. But series such as Av(12453) are not as clear cut. So if there is another singularity type, or even stretched exponentials
with additional multiplicative logarithms corrections, we are not testing for that.

\section{Program to generate coefficients.}

We first wrote a general-purpose algorithm which produced data for all
Wilf classes for patterns up to O$(x^{16})$  (which means 17 terms in the
generating function). Our algorithm takes time $O(T\cdot s_n)$, where
$T$ is the time taken to check whether a pattern is avoided, and $s_n$
is the number of PAPs of length $n.$ Our algorithm has essentially
zero memory requirements, while the testing time $T$ is $O(n^{k-1})$
for a pattern of length $k$. 

We then
learnt of Kuszmaul's algorithm \cite{K17} which also
takes time O$(T\cdot s_n),$ but $T$ is now O$(k),$ though the trade-off
is significantly greater memory requirements. However, for the
permutations up to length 16 the memory requirements are still modest,
and Kuszmaul's algorithm is faster than ours -- typically by a factor
of 10 -- so we see no point in describing ours.

Inoue and Minao proposed a radically different algorithm \cite{IM14} which does not examine each permutation
individually. It uses a permutation set description technique
called $\pi$DD which has established reasonably efficient algorithms for
various set operations on such representations \cite{ITM13}. We describe this in the following subsection.

\subsection{The Rot-$\pi$DD algorithm}

Imagine one has variables $x_1$, $x_2$, $\dots x_n$.

A {\it Binary Decision Diagram} (BDD) \cite{B86}
is an often compact method of representing
sets of sets of these variables. It is a tree representation where each
node is one of the variables. To determine whether a
particular
set $S\subseteq\{x_1, x_2, \dots x_n\}$ is in the represented set of
sets,
start at the root node, and traverse
down the tree.
 At each node  $x_i$ it will
have two children. Take the left if $x_i \notin S$ and the right if
$x_i \in S$.
Eventually one will reach a leaf which is $0$ or $1$ representing
respectively
that $S$ is not or is in the represented set of sets. As a complete
binary tree,
there are $2^n$ leaves and this is not compact; however merging of
isomorphic
subgraphs and eliminating nodes whose children are isomorphic can often
reduce
the representation to a more efficient version.

Importantly, there are well known reasonably efficient
algorithms for various operations on a
BDD, including set intersection, union, and cardinality.

There is a slight improvement on a BDD called a {\it Zero Suppressed
Decision Diagram} (ZDD) \cite{M93}
which can be more compact by suppressing nodes whose positive edge
points
to the zero leaf. More precisely, the compression ratios of BDD and ZDD depend on the particular set of sets they represent. 
BDD was historically invented to represent a Boolean function and has a better compression ratio on Boolean functions empirically. 
On the other hand, ZDD can compress a set of sparse sets well empirically. Since permutation sets on combinatorial problems have some 
mathematical structure and tend to be represented as sparse sets of operations, we selected ZDD as the base structure of $\pi$DD rather than BDD.

A permutation can be represented as a set of operations that
generate
the permutation. This means that a set of permutations can be
represented
as a BDD or better ZDD where each variable is one of the operations
that
generate the permutation. This is called
a {\it Permutation Decision Diagram} ($\pi$DD). Again the point is that
the representation of the $\pi$DD may be significantly more compact
than a list of the permutations it represents \cite{M11}.
The set of operations that are used in the standard representation of a $\pi$DD
are element exchanges.
Any $n$ element permutation can be written as the
composition
of up to $n-1$ two element exchanges.  A Rot-$\pi$DD   \cite{IM14} is the same idea
except that
the variables are now, instead of a swap of elements $i$ and $j$, a
rotation
of the elements between $i$ and $j$ inclusive. Conceptually this
is a very similar approach; in practice however the constructed
sets for many problems using the Rot-$\pi$DD representation is
significantly more compact than the $\pi$DD representation.
See \cite{I21} for a detailed
description of the algorithm.

The 
publications \cite{ITM13} and \cite{I17} describe a method of constructing all permutations that do
{\it not} avoid a pattern by constructing representations of
permutations
that  include the pattern in each possible specific choice of
elements in the permutation. These sets turn out to be significantly
more compact with the Rot-$\pi$DD representation than the $\pi$DD
representation.  The union of these sets (performed with standard
ZDD algorithms) is the set of all non-PAPs. Subtracting from the
set of all permutations then gives the set of all PAPs. 

\subsection{Comparison with previous algorithms}

The $\pi$DD algorithm has significantly different computational
requirements to the $O(T\cdot s_n)$ algorithms which visit each pattern
avoiding permutation, and can in principle be significantly faster than
any such algorithm, at the cost of a higher memory consumption.

While good theoretical complexity analysis is not available for the $\pi$DD
or Rot-$\pi$DD algorithms, given current computer speeds and memory
capacities, the Rot-$\pi$DD algorithm is multiple orders of magnitude
faster than algorithms like Kuszmaul's or our prior algorithm that
inspect each pattern avoiding permutation individually, typically
producing an extra 8 or 10 terms. To put this into context, previous algorithms would
require years of CPU time to produce series of this length, rather than hours required by the Rot-$\pi$DD algorithm.

The data in this paper were computed using Inoue's Rot-$\pi$DD program, and the
early terms were independently checked by both our early simple
algorithm and Kuszmaul's program. The first 20 digits of the later terms
were confirmed by the method of series extension, described below. Full details can be found in \cite{I21}.

The calculations were done using the Research Computing Services facilities hosted at the University of Melbourne, more precisely the Melbourne Research Cloud,  based on the OpenStack cloud computing platform. The virtual machine
used had 2TB of memory and 48 cores and was running on a physical node with an AMD EPYC 7702 64-Core Processor.

\section{Stieltjes moment sequences}
\label{stieltjesn}
The classical {\em Stieltjes moment problem} considers a numerical sequence ${\bf a} \equiv \{a_n\}, \, n \ge 0$ in which $a_n$ can be expressed as the integral$$a_n = \int_\Gamma x^n d\rho(x) $$ for all $n \ge 0,$ where the support $\Gamma \subseteq {\mathbb R},$ and $\rho$ is a measure. If the measure $\rho$ is differentiable, then it has a {\em density} or {\em density function} $\mu(x)=\rho'(x).$  In which case the above equation becomes
$$a_n = \int_\Gamma x^n \mu(x) dx. $$ 
There are several equivalent conditions that the sequence ${\bf a}$ must satisfy in order to be a Stieltjes moment sequence, or equivalently, for such a density function to exist, which must of course be non-negative. 

The Hankel matrix $H_n^\infty({\bf a})$ is defined as
\[
H_n^\infty({\bf a})=
  \begin{bmatrix}
    a_n & a_{n+1} & a_{n+2} & \ldots \\
    a_{n+1} & a_{n+2}& a_{n+3} & \ldots\\
   a_{n+2} & a_{n+3}& a_{n+4} & \ldots\\
\vdots & \vdots & \vdots & \ddots \\
  \end{bmatrix}
\]
The following theorem was proved in part by Stieltjes and in part by Gantmakher and Krein. In particular, the properties (a) and (d) (below) were shown to be equivalent by Stieltjes \cite{St}, while these were later shown to be equivalent to (b) and (c) by Gantmakher and Krein\cite{GK}. 
{\theorem For a sequence ${\bf a} \equiv \{a_n\}, \, n \ge 0$, the following are equivalent:
(a) There exists a positive measure $\rho$ on $\Gamma \in [0,\infty)$ such that
$$a_n = \int_\Gamma x^n d\rho(x). $$
(b) The matrices $H_0^\infty({\bf a})$ and $H_1^\infty({\bf a})$ are both positive semidefinite.\\
(c) The matrix $H_0^\infty({\bf a})$ is totally positive (all of its minors are non-negative). \\
(d) There exists a sequence of real numbers $\alpha_0, \, \alpha_1,\ldots \ge 0$  such that the generating function $A(x)$ for the sequence $a_0, a_1, \ldots$  satisfies

\[A(x) = \sum_{n=0}^\infty a_n x^n =  \cfrac{\alpha_0}{1 -\cfrac{\alpha_1 x}{1- \cfrac{\alpha_2x}{
      \begin{array}{@{}c@{}c@{}c@{}}
        1 - \ldots 
      \end{array}
    }}} \]}

 A sequence that satisfies the equivalent conditions of Theorem 1 is called a {\em Stieltjes moment sequence}.

One reason for attempting to identify combinatorial sequences as Stieltjes moment sequences is that such sequences are log-convex.
To see that a Stieltjes moment sequence is log-convex, it suffices to observe that for each $n \ge 1$, the expression $a_{n+1}a_{n-1}-a_n^2 \ge 0$ is a minor of $H_0^\infty({\bf a})$, so this expression is non-negative by condition (c).


 Log-convexity of the sequence ${\bf a}$ implies that the ratios $\frac {a_n}{a_{n-1}}$ provide lower bounds to the growth constant $\mu$ of the sequence. In the case that ${\bf a}$ is a Stieltjes moment sequence, we can use the above properties to compute stronger lower bounds for $\mu,$ using a method first given by Haagerup, Haagerup and Ramirez-Solano \cite{HHR15}.

Using the coefficients $a_0, a_1, \ldots a_n,$  we calculate the terms $\alpha_0, \alpha_1, \ldots \alpha_n$ in the continued fraction representation above. It is easy to see that the coefficients of $A(x)$ are non-decreasing in each $\alpha_j.$ Hence $A(x)$ is (coefficient-wise) bounded below by the generating function $A_n(x)$, defined by setting $\alpha_n, \alpha_{n+1},  \alpha_{n+2},\ldots $ to 0. Therefore, the growth rate $\mu_n$ of $A_n(x)$ is no greater than the growth rate $\mu$ of $A(x).$ The growth rates $\mu_1, \mu_2, \ldots $ clearly form a non-decreasing sequence, and, since the coefficients of $A_n(x)$ are log-convex, $\mu_n \ge a_n/a_{n-1}$. It follows that this sequence $\mu_1, \mu_2, \ldots $ of lower bounds converges to the exponential growth rate $\mu$ of {\bf a}.

If we assume further that the sequences $\alpha_0, \alpha_2, \alpha_4 \ldots $ and $\alpha_1, \alpha_3, \alpha_5 \ldots $ are non-decreasing, as we find empirically in all of the cases we consider, we can get stronger lower bounds for the growth rate by setting $\alpha_{n+1}, \alpha_{n+3}, \ldots$  to  $\alpha_{n-1}$ and $\alpha_{n+2}, \alpha_{n+4}, \ldots$    to $\alpha_{n}$. For this sequence the exponential growth rate of the corresponding sequence ${\bf a}$ is  $b_n \equiv (\sqrt{\alpha_n}+\sqrt{\alpha_{n-1}})^2.$ By the method with which we constructed this bound, it is clear that
$(\sqrt{\alpha_n}+\sqrt{\alpha_{n-1}})^2 \ge \mu_n.$  Hence, the lower bounds $(\sqrt{\alpha_n}+\sqrt{\alpha_{n-1}})^2$ converge to the growth constant $\mu.$ In all cases considered here these bounds are monotone in $n,$ and can be numerically extrapolated. In particular, if the ratios of the original coefficients, $$r_n =\frac{a_n}{a_{n-1}} \sim \mu \left (1+\frac{c}{n^\theta} + o \left ( \frac{1}{n^\theta} \right ) \right ), $$ then we have heuristic arguments, (not yet a proof), that $b_n \sim const. \cdot n^\beta,$ where $\beta=2\theta/(2-\theta).$  This gives us the appropriate power against which to extrapolate the bounds. Also, if $\alpha_n \le \alpha_{n+2}$ for each $n$ and the limit
$$\lim_{n \to \infty} \alpha_n$$ exists, then it is equal to $\mu/4.$ We have calculated these ``bounds" for all sixteen Wilf classes.

\section{Classical patterns of length 5}

The 120 possible classical PAPs of length five belong to just 16 Wilf classes. Only one of these 16 cases is solved, that of $Av(12345),$ in
the sense that the generating function can be  explicitly written down.  In eqn. (34) in \cite{BEGM20} the solution is given as a second-order linear  inhomogeneous ODE,
and furthermore the generating function is shown to be a Stieltjes moment sequence. Furthermore, it is simply related to the moments of a 4-step random walk in two dimensions, as
well as to the number of $2n$-step polygons on the diamond lattice.

 In Tables \ref{tab:terms1} and \ref{tab:terms2} we give the coefficients to various orders from 23 to 27 for all 16 Wilf classes, as generated by our program.

We have  analysed the series in order to estimate the growth constants and the associated exponent of the sub-dominant term.
We also identified those cases which we believe have a stretched-exponential singularity, as observed in the case of $Av(1324)$ PAPs, and discussed above.

We have also identified the permutation classes of all 120 possible length-5 classical permutations. The obvious equality
$$s_n(n_1n_2n_3n_4n_5)=s_n(n_5n_4n_3n_2n_1)$$ reduces the identification to 60 possible classes. Further reductions follow from
the results \cite{BW00} that $$s_n(54n_3n_4n_5)=s_n(45n_3n_4n_5)$$ and $$s_n(543n_4n_5)=s_n(345n_4n_5),$$  the result proved in \cite{SW02}
that $$s_n(354n_4n_5)=s_n(435n_4n_5),$$ the result \cite{BWX07} that $$s_n(12\cdots k \tau)=s_n(k k-1\cdots 1 \tau),$$ where $\tau \in [k+1..l],$ and 
$$s_n(\tau k k-1)=s_n(\tau k-1 k),$$ proved in \cite{R03}.
Further identification of Wilf classes using the Wilf-equivalence between two-layer permutations and monotone permutations of the same length, as discussed in \cite{V15}, has also been made.
 The remaining equivalences were established by direct enumeration of the outstanding classes
using our computer program. The results are shown in table \ref{tab:classes} below. Anticipating our numerical results, we have listed the 16 classes in increasing order of the (estimated) growth constants. The result of this is that they are {\em not} necessarily listed in increasing order of the last-known coefficient. Of course, if sufficiently many coefficients were known, these two orderings must be identical.

\begin{table}[htbp] 
   \begin{center}
   \small
   \begin{tabular}{|r r r r r r r r r r r|} \hline
  25314&41352&&&&&&&&&\\\hline
31524&42513&24153&35142&&&&&&&\\\hline
35214&41253&23514&41532&25134&43152&25413&31452&&&\\\hline
43251&15234&13452&25431&23415&51432&41235&53214&&&\\\hline
34215&51243&14532&23541&15423&32451&43125&52134&&&\\\hline
53124&42135&13542&24531&15243&34251&32415&51423&&&\\\hline
32541&14523&34125&52143&&&&&&&\\\hline
35124&42153&24513&31542&25143&34152&41523&32514&&&\\\hline
31245&54213&12453&35421&12534&43521&21453&35412&21534&43512&\\
23145&54132&23154&45132&31254&45213&&&&&\\\hline
42351&15324&14352&25341&24315&51342&41325&52314&&&\\\hline
42315&51324&15342&24351&&&&&&&\\\hline
12345&54321&45321&12354&12543&34521&21345&54312&21354&45312&\\
21543&34512&23451&15432&32145&54123&32154&45123&43215&51234&\\\hline
35241&14253&13524&42531&24135&53142&31425&52413&&&\\\hline
53241&14235&13425&52431&&&&&&&\\\hline
53421&12435&21435&53412&13245&54231&13254&45231&&&\\\hline
52341&14325&&&&&&&&&\\\hline

\hline
      \end{tabular}
        \caption{The distribution of the 120 possible length-5 pattern-avoiding permutations among the 16 Wilf classes. Note that the entries for the classes $Av(31245)$ and $Av(12345)$ are split across two lines of the table.}
   \label{tab:classes}
   \end{center}
\end{table}

\begin{table}[htbp]
\begin{center}
\tiny
\begin{tabular}{rrrrr} \hline
$n$& Av(25314) A256195&      Av(31524) A256196&    Av(35214) A256197&      Av(43251) A256203\\
7  & 4578&                   4579&                  4579&                   4581\\
8  & 33184&                  33216&                 33218&                  33283\\
9  & 258757&                 259401&                259483&                 260805\\
10 & 2136978&                2147525&               2149558&                2171393\\
11 & 18478134&               18632512&              18672277&               18994464\\
12 & 165857600&              167969934&             168648090&              173094540\\
13 & 1535336290&             1563027614&            1573625606&             1632480259\\
14 & 14584260700&            14937175825&           15093309024&            15851668551\\
15 & 141603589300&           146016423713&          148223240022&           157824649955\\
16 & 1400942032152&          1455402205257&         1485673163882&          1605839173312\\
17 & 14087464765300&         14753501614541&        15159644212775&         16652321922596\\
18 & 143689133196008&        151783381341695&       157142812302992&        175596537163347\\
19 & 1484090443264936&       1582029822426003&      1651865171372967&       1879357191026029\\
20 & 15499968503875136&      16681492660789425&     17582693993265148&      20382942631855557\\
21 & 163501005435759505&     177726496203056670&    189269329080075275&     223719376672365073\\
22 & 1740170514634463426&    1911230701872865231&   2058215511081891400&    2482094083780961295\\
23 & 18671118911254798454&   20726637978574528119&  22589841589522026553&   27808544385768051233\\
24 & 201805434191401310152&  226497541235099049284& 250032335770049925668&  314346011323933283258\\
25 & 2195829593847519231848& 2492440846906157577367&2788899325208909923567& 3582440933577530273836 \\
26 & 24039330044242839545400&&                      31329479505464363566868&41134198972534449502215\\
27 & &                       &                      &                       475581766378016525358137 \\
   & &&&\\
$n$& Av(34215) A256205&       Av(53124) A256199 &    Av(32541) A256204 &     Av(35124) A256198 \\
7  & 4581&                    4580&                  4581&                   4580\\
8  & 33285&                   33252&                 33284&                  33249\\
9  & 260886&                  260202&                260847&                 260092\\
10 & 2173374&                 2161837&               2172454&                2159381\\
11 & 19032746&                18858720&              19015582&               18815124\\
12 & 173741467&               171285237&             173461305&              170605392\\
13 & 1642533692&              1609282391&            1638327423&             1599499163\\
14 & 15999488304&             15561356705&           15939733122&            15427796984\\
15 & 159917206735&            154246419725&          159099927785&           152487271455\\
16 & 1634681988983&           1562151687940&         1623799173782&          1539554179950\\
17 & 17042352950764&          16121960812335&        16900201391546&         15836801521762\\
18 & 180798150762914&         169178376076607&       178967276844263&        165625811815111\\
19 & 1948027746498015&        1801800479418116&      1924689980696921&       1757953168747511\\
20 & 21282786390947602&       19446010522240384&     20987593594256974&      18908510233855411\\ 
21 & 235446451502773103&      212394673429250090&    231734179050033660&     205838673911323648\\
22 & 2634317655935012208&     2345064355131025130&   2587835777992844938 &   2265393020812413370\\
23 & 29778833170013213300&    26148064110299271293&  29198736751160012102&   25182471016157568626\\
24 & 339796984870771392635 &  294190661855648481179 &332575357468837097628&  282511039355447739772\\
25 & 3910755764784092153311&  3337335970674441425688&3821024002600674745994& 3196265588333257586119\\
26 & 45365839293606522375359& &                      44252507544177176282956&36445643066828928379492\\
27 & 530098601158553050947014&&                      &                       418600631627670270370879\\ 
\hline
\end{tabular}
\caption{Coefficients from order 6 to 27 for the half of the 16 Wilf classes
with smaller growth constant $\mu$. 
The initial terms in all sequences are $1,1,2,6,24,119,694.$ 
The classes are arranged in increasing order with respect to $\mu$,
reading left to right from the top row.}
\label{tab:terms1}
\end{center}
\end{table}

\begin{table}[htbp]
\begin{center}
\tiny
\begin{tabular}{rrrrr} \hline
$n$ & Av(31245) A116485 &      Av(42351) A256200 &      Av(42315) A256206 &     Av(12345) A047889 \\
7  & 4581&                    4580&                    4581&                   4582\\
8  & 33286&                   33252&                   33287&                  33324\\
9  & 260927&                  260204&                  260967&                 261808\\
10 & 2174398&                 2161930&                 2175379&                2190688\\
11 & 19053058&                18861307&                19072271&               19318688\\
12 & 174094868&               171341565&               174426353&              178108704\\
13 & 1648198050&              1610345257&              1653484169&             1705985883\\
14 & 16085475576&             15579644765&             16165513608&            16891621166\\
15 & 161174636600&            154541844196&            162344264849&           172188608886\\
16 & 1652590573612&           1566713947713&           1669261805697&          1801013405436\\
17 & 17292601075489&          16190122718865&          17526017429722&         19274897768196\\
18 & 184246699159418&         170171678529883&         187472773174466&        210573149141896\\
19 & 1995064785620557&        1816001425551270&        2039233971499931&       2343553478425816\\
20 & 21919480341617102&       19646035298044543&       22520066337196663&      26525044132374656\\
21 & 244015986016996763&      215179180467834605&      252141732452056894 &    304856947930144656\\
22 & 2749174129340156922&     2383465957654163227&     2858721279079666465&    3553266124166899872\\
23 & 31313478171012371344&    26673704385975326866&    32786666580814894741&   41952101272633801376\\
24 & 360255986786421416732&   301342110309622207830&   380034587229949049485&  501228159413699278144\\
25 & 4183070452633759090955&  3434155564505269412223&  4448342812221497172384& 6054582181256780696704\\
26 & 48986523769015357032198& 39453283522954708152659& 52542550112506952504622&73884542290182291304704\\
27 & 578206680078321677926243&456668245606432017686247&&                       910193895170720544149248\\
   & &&&\\
$n$& Av(35241) A256201 &     Av(53241) A256202 &      Av(53421) A256207 &     Av(52341) A256208\\
7  & 4580&                   4580&                    4582&                   4582\\
8  & 33254&                  33256&                   33325&                  33325\\
9  & 260285&                 260370&                  261853&                 261863\\
10 & 2163930&                2166120&                 2191902&                2192390\\
11 & 18900534&               18945144&                19344408&               19358590\\
12 & 172016256&              172810050&               178582940&              178904675\\
13 & 1621031261&             1633997788&              1713999264&             1720317763\\
14 & 15739870457&            15939893003&             17019444969&            17132629082\\
15 & 156855197297&           159820729208&            174149184184&           176055309619\\
16 & 1599233708733&          1641980432159&           1830279810276&          1861037944163\\
17 & 16638560125635&         17242378256155&          19703572779755&         20185165186517\\
18 & 176269571712376&        184674461615836&         216769635980879&        224150069984572\\
19 & 1898076560618372&       2013829450204384&        2432308876304981&       2543698932578158\\
20 & 20742488003444465&      22324460502429244&       27788506478197951&      29451619807433107\\
21 & 229747253093647567 &    251250502143635615&      322770995262901091&     347417296695040510\\
22 & 2576270755655436479&    2867467023751687892&     3806657237502632706 &   4170088041714300134\\
23 & 29218474225923168362&   33152272498223444540&    45532086120583546634&   50874753262007210667\\
24 & 334868638387692996919 & 387935538721724466875&   551794232925251495478&  \\
25 & 3875365114838257507148& 4590792008759551665335&  6769119579399164598190 &\\
26 & 45256353903547788096108&54901471673327772683658 &83991144346393508063125&\\
27 & &                       &                        &                       \\
\hline
\end{tabular}
\caption{Coefficients from order 6 to 27 for the half of the 16 Wilf classes
with larger growth constant $\mu$. 
The initial terms in all sequences are $1,1,2,6,24,119,694.$ 
The classes are arranged in increasing order with respect to $\mu$,
reading left to right from the top row.}
\label{tab:terms2}
\end{center}
\end{table}

\section{Series analysis}
The method of series analysis  has, for many years, been a powerful tool in the study of a variety of problems in statistical mechanics, combinatorics, fluid mechanics and computer science. In essence, the problem is the following: Given the first $N$ coefficients of the series expansion of some function, (where $N$ is typically as low as 5 or 6, or as high as 100,000 or more), determine the asymptotic form of the coefficients, subject to some underlying assumption about the asymptotic form, or, equivalently, the nature of the singularity of the function.

A typical example is the generating function of self-avoiding walks (SAWs) in dimension two or three. This is believed to behave as 
\begin{equation}\label{generic}
F(z) = \sum_n c_n z^n \sim C \cdot (1 - z/z_c)^{-\gamma}.
\end{equation}
  In this case, among regular two-dimensional lattices, the value of $z_c$ is only known for the hexagonal lattice \cite{DC10}, while $\gamma=43/32$ is believed to be the correct exponent value for all two-dimensional lattices, but this has not been proved. 

The method of series analysis is used when one or more of the critical parameters is not known. For example, for the three-dimensional versions of the above problems, none of the quantities $C,$ $z_c$ or $\gamma$ are exactly known.
 From the binomial theorem it follows from (\ref{generic}) that 
\begin{equation}\label{asymp1}
c_n \sim \frac{C }{\Gamma(\gamma)}\cdot  z_c^{-n} \cdot  n^{\gamma-1}
\end{equation}
 Here $C, \,\, z_c, \,\, {\rm and} \,\, \gamma$ are referred to as the critical amplitude, the critical point (usually the radius of convergence) and the critical exponent respectively. In combinatorics one often refers to the {\em growth constant $\mu$}, which is $\mu=1/z_c,$ as the coefficients are dominated by the term $\mu^n.$

Obtaining these coefficients is typically a problem of exponential complexity, as is the case with our algorithm, described in Sec. 2. The usual consequence is that that fewer than 50 terms are known (and in some cases far fewer).

\section{Ratio Method}
The ratio method was perhaps the earliest systematic
method of series analysis employed,
and is still the most useful method when only a small number of terms are known.
From equation (\ref{asymp1}),
it follows that the {\it ratio} of successive terms
\begin{equation} \label{ratios}
r_n = \frac{c_n}{c_{n-1}}=\frac{1}{z_c}\left (1 + \frac{\gamma -1}{n} + {\rm o}(\frac{1}{n})\right ).
\end{equation}
 It is then natural to plot the successive ratios $r_n$ against $1/n.$
If the correction terms ${\rm o}(\frac{1}{n})$ can be ignored\footnote{For a purely algebraic singularity eqn. (\ref{generic}), with no confluent terms, the correction term will be ${\rm O}(\frac{1}{n^2}).$}, such a plot will be linear,
with gradient $\frac{\gamma-1}{z_c},$ and intercept $\mu=1/z_c$ at $1/n = 0.$

\begin{figure}[h!] 
\begin{minipage}[t]{0.45\textwidth} 
\centerline{
\includegraphics[width=\textwidth]{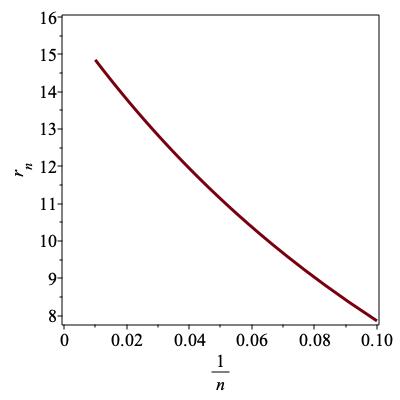}
}
\caption{
Plot of ratios against $1/n$ for $Av(12345)$ PAPs. A straight line through the
last few data points intercepts the ratios axis at $\mu=1/z_c.$
} 
\label{fig:ratt}
\end{minipage}
\hspace{0.05\textwidth}
\begin{minipage}[t]{0.45\textwidth} 
\centerline {\includegraphics[width=\textwidth]{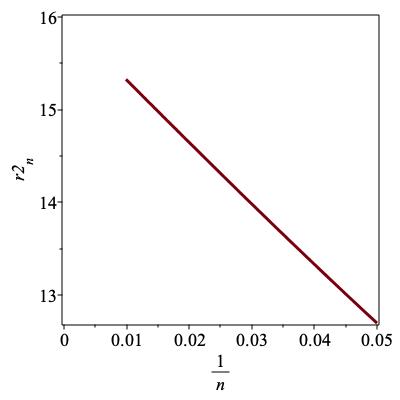}}
\caption{Plot of modified ratios $r2_n$ against $1/n$ for $Av(12345)$ PAPs, showing linearity.} \label{fig:modrat}
\end{minipage}
\end{figure}

As an example, we apply the ratio method to the solved case of $Av(12345)$ PAPs \cite{BEGM20}. We will use 101 coefficients, to align with our analysis of other PAPs subsequently.
Plotting successive ratios against $1/n$ results in the plot shown in Fig.~\ref{fig:ratt}.
The critical point is known to be at $z_c  = 1/16.$

\begin{figure}[h!] 
\begin{minipage}[t]{0.45\textwidth} 
\centerline{
\includegraphics[width=\textwidth]{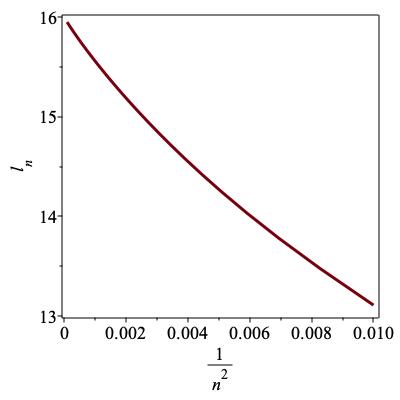}
}
\caption{
Plot of linear intercepts $l_n$ against $1/n^2$ for $Av(12345)$ PAPs. A straight line through the
last few data points intercepts the vertical axis at $\mu=1/z_c.$
} \label{fig:ln}

\end{minipage}
\hspace{0.05\textwidth}
\begin{minipage}[t]{0.45\textwidth} 
\centerline {\includegraphics[width=\textwidth]{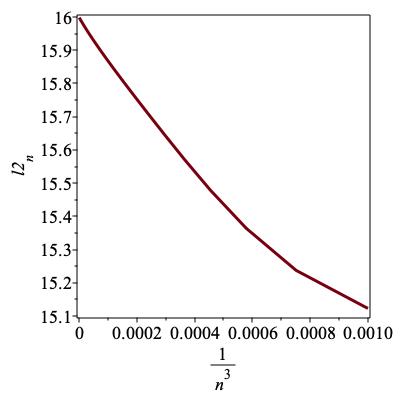}}
\caption{
Plot of quadratic intercepts $l2_n$ against $1/n^3$ for $Av(12345)$ PAPs. A straight line through the
last few data points intercepts the vertical axis at $\mu=1/z_c.$
} \label{fig:ln2}
\end{minipage}
\end{figure}

 From the figure one sees that the locus of points still displays some curvature. There are two possible explanations. One is that there is a stretched exponential term, the other is that we have a pure power law, but that the term of order $1/n^2$ has a coefficient much greater than the coefficient of the term of order $1/n,$ so makes a strong contribution for small values of $n$. To see if this is in fact the case, we can eliminate the term of order $1/n^2$ by forming the modified ratios:
 \BE \label{eqn:modrat}
 r2_n \equiv \frac{n^2 r_n - (n-1)^2 r_{n-1}}{2n} = \mu \left ( 1 + \frac{g-1}{2n} + O \left ( \frac{1}{n^3} \right ) \right ).
 \EE
 We show in Fig.~\ref{fig:modrat} a plot of $r2_n$ against $1/n,$ which now appears to be completely linear, supporting the second explanation for the initial curvature of the ratio plot. Indeed, further numerical investigation of the ratios show that they behave as
 $$r_n = \mu \left ( 1 - \frac{7.5}{2n} +  \frac{131.25}{n^2} + O \left ( \frac{1}{n^3} \right ) \right ),$$ showing the very large O$(1/n^2)$ term compared to the order $1/n$ term.

\begin{figure}[h!] 
\begin{minipage}[t]{0.45\textwidth} 
\centerline {
\includegraphics[width=\textwidth]{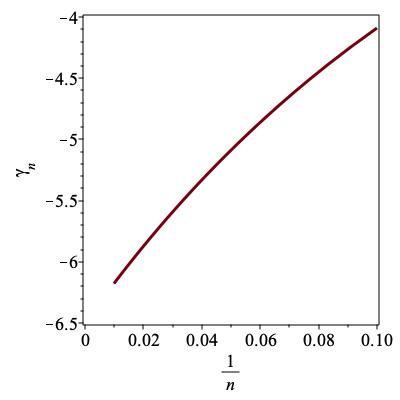}
}
\caption{Plot of exponent estimates $\gamma_n$ against $1/n$ for $Av(12345)$ PAPs, assuming $\mu=16.$ Visual extrapolation suggests $\gamma \approx -6.5.$}
 \label{fig:gamn}
\end{minipage}
\hspace{0.05\textwidth}
\begin{minipage}[t]{0.45\textwidth} 
\centerline {
\includegraphics[width=\textwidth]{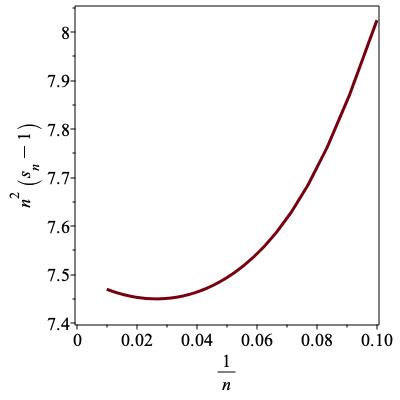}}
\caption{Plot of  $n^2(s_n-1)$ against $1/n$ for $Av(12345)$ PAPs, which should go to $-g=7.5$ as $n \to \infty.$}
 \label{fig:A12345div}
\end{minipage}
\end{figure}

Visual extrapolation to $1/z_c \approx 16$ is quite obvious. A straight line drawn through the
last $4-6$ data points intercepts the horizontal axis around $1/n \approx 0.07.$ Thus
the gradient is approximately $\frac{16-7.8}{-0.07} \approx -117,$
from which we conclude that the exponent $\gamma - 1 \approx -117 \cdot z_c \approx -7.3.$ It is known
that the exact value is $\gamma-1 = -7.5,$ which is in approximate agreement with this
simple graphical analysis.

Linear intercepts $l_n$ eliminate the $O\left ( \frac{1}{n} \right )$ term in eqn. (\ref{ratios}), so in the case of a pure power-law singularity, one has
$$l_n \equiv nr_n - (n-1)r_{n-1} = \mu \left (1+ \frac{c}{n^2} + O\left (\frac{1}{n^3} \right ) \right ).$$ The linear intercepts in this case are shown in Fig.~\ref{fig:ln}, and it is clear that
they are converging to the known limit $\mu=16$ rather more rapidly than are the plain ratios. In the case of a pure power-law, where sub-dominant terms in the ratios decrease by 
successive factors of $1/n,$ this process can be continued. For example, we can eliminate terms of order $1/n^2$ by forming quadratic estimators
$$l2_n \equiv \frac{n^2l_n - (n-1)^2 l_{n-1}}{2l-1} = \mu \left (1+ \frac{c}{n^3} + O\left (\frac{1}{n^4} \right ) \right ).$$ The quadratic intercepts are shown in Fig.~\ref{fig:ln2}, and it is clear that
they are converging to the known limit $\mu=16$ even more rapidly than are the linear intercepts or ratios. 

Various refinements of the method can be readily derived. If the critical point
is known exactly, it follows from equation (\ref{ratios}) that estimators of the exponent
$\gamma$ are given by
$$ \gamma_n=n(z_c\cdot r_n-1)+1 = \gamma +  {\rm o}(1).$$

The estimators $\gamma_n$ are shown in Fig.~\ref{fig:gamn}, and it can be seen that they are plausibly going to a limit, as $1/n \to 0,$ of -6.5, which is the exact value
of the exponent in this case.

If the critical point is not known exactly, one can still estimate the exponent $\gamma.$ From eqn (\ref{ratios}) it follows that 
\BE \label{eq:exp}
\delta_n=1+n^2\left ( 1-\frac{r_n}{r_{n-1}} \right )= \gamma +  {\rm o}(1).
\EE

\begin{figure}[h!] 
\begin{minipage}[t]{0.45\textwidth} 
\centerline {\includegraphics[width=\textwidth]{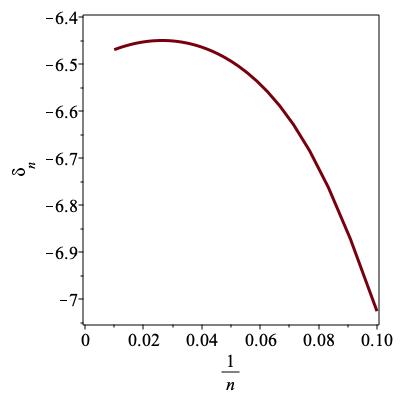}}
\caption{
Plot of exponent estimates $\delta_n$ against $1/n$ for $Av(12345)$ PAPs.
} \label{fig:deln}
\end{minipage}
\hspace{0.05\textwidth}
\begin{minipage}[t]{0.45\textwidth} 
\centerline {
\includegraphics[width=\textwidth]{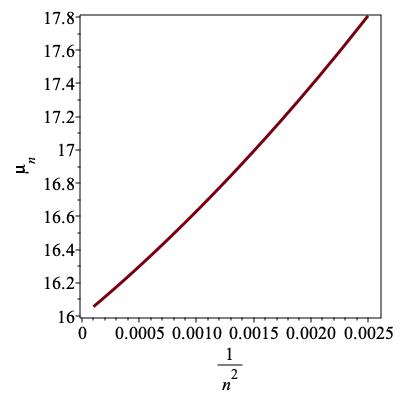}
}
\caption{Plot of growth-constant estimates $\mu_n$ against $1/n^2$ for $Av(12345)$ PAPs.}
 \label{fig:mun}
\end{minipage}
\end{figure}

The estimators $\delta_n$ are shown in Fig.~\ref{fig:deln}, and it can be seen that they too are plausibly going to a limit, as $1/n \to 0,$ of $\gamma=-6.5.$

Similarly, if the exponent $\gamma$ is known, estimators of the growth constant  $\mu$
are given by  $$\mu_n =  \frac{n r_n}{n+\gamma-1}= \mu + {\rm o}(1/n).$$

The estimators $\mu_n$ are shown in Fig.~\ref{fig:mun}, and it can be seen that they  are very plausibly going to a limit, as $1/n \to 0,$ of $\mu = 16.$

The explanation of Fig.~\ref{fig:A12345div} is given in Section \ref{sec:se} after eqn. (\ref{eqn:div}). The figure  provides compelling evidence for a pure power-law singularity in this case.

\section{Functions with non-power-law singularities.}\label{sec:nonalg}

A number of solved, and, we believe, unsolved problems that arise in lattice critical phenomena and algebraic combinatorics have coefficients with a more complex asymptotic form, with a sub-dominant  term ${\rm O}(\mu_1^{n^\sigma})$ as well as a power-law term ${\rm O}(n^g).$  Perhaps the best-known example of this sort of behaviour is the number of partitions of the integers -- though in that case the leading exponential growth term $\mu^n$ is absent (or equivalently $\mu=1$).
The  form of the coefficients $c_n$ in the general case is
\BE \label{stretch}
c_n \sim C \cdot \mu^n \cdot \mu_1^{n^\sigma} \cdot n^g.
\EE

An example from combinatorics is given by Dyck paths enumerated not just by
length, but also by height (defined to be the maximum vertical distance of the
path from the horizontal axis). Let $d_{n,h}$ be the number of Dyck paths of
length $2n$ and height $h.$ The OGF is then\footnote{One of us (AJG) posed this problem at an Oberwolfach meeting in March 2014. Within 24 hours Brendan McKay produced this solution. See also \cite{NP13}.}
\BE
D(x,y) = \sum_{n,h} d_{n,h}x^{2n}y^h,\,\,{\rm and} \,\, [x^{2n}]D(x,y)=\sum_{h=1}^n d_{n,h}y^h.
\EE
For $y < 1$
let $A=2^{5/3}\pi^{5/6}/\sqrt{3},$ $E=3\left ( \frac{\pi}{2} \right )^{2/3}$ and $r=-\log{y}.$
Then one finds that $[x^{2n}]D(x,y)$ is given by eqn. (\ref{stretch}) with $C=\frac{1-y}{y^2}r^{1/3}A,$ $\mu=4,$ $\mu_1=\exp(-Er^{2/3}),$ $ \sigma=1/3,$ and $g=-5/6.$

Note however that such singular behaviour {\em can} arise from D-finite ODEs as the asymptotic form of the coefficients of some generating function at an {\em irregular} singular point. The extraction of the asymptotic behaviour from an ODE at an irregular singular point is not something that has been automated, as each case must be treated individually, and usually involves some subtleties. 

But the general solutions invariably include a factor 
\BE \label{irregular}
\mu^n \cdot \exp(R(n^{1/\rho}))\cdot n^g,
\EE
 where $\rho$ is a positive integer and $R$ is a polynomial

Applying the ratio method to such singularities requires some significant changes. These were first developed in \cite{G15}, where further details and more examples can be found. 
In the next subsection we give a summary, including as much detail as is needed for our analysis.

\subsection{Ratio method for stretched-exponential singularities.} \label{sec:se}

If 

\BE \label{eq:an}
b_n \sim C \cdot \mu^n \cdot \mu_1^{n^\sigma} \cdot n^g,
\EE
 then the ratio of successive coefficients $r_n = b_n/b_{n-1},$ is
\begin{multline} \label{eq:rn}
r_n = \mu \left (1 + \frac{\sigma \log \mu_1}{n^{1-\sigma}} + \frac{g}{n} + \frac{\sigma^2 \log^2 \mu_1}{2n^{2-2\sigma}} + \frac {(\sigma-\sigma^2)\log \mu_1+2g\sigma \log \mu_1}{2n^{2-\sigma}} \right . \\
 \left . {}+ \frac{\sigma^3 \log^3 \mu_1}{6n^{3-3\sigma}} +{\rm O}(n^{2\sigma-3}) + {\rm O}(n^{-2}) \right ).
\end{multline}

It is usually the case that $\sigma$ takes the simple values $1/2,$ $1/3,$ $1/4$ etc.\footnote{In statistical mechanical models, the value of the exponent $\sigma$ is simply related to the fractal dimension $d_f$  of the object through $\sigma = 1/(1+d_f).$ }. Recall too, from eqn. (\ref{irregular}) that if these asymptotics arise as an irregular singular point of a D-finite ODE, the exponent must be of the form $1/\rho,$ where $\rho$ is a positive integer \cite{FS09}.

When $\sigma = \frac{1}{2},$ (\ref{eq:rn}) specialises to
\BE \label{eq:half}
r_n = \mu \left (1 + \frac{ \log \mu_1}{2\sqrt{n}} + \frac{g+\frac{1}{8}\log^2 \mu_1}{n} + \frac{\log^3\mu_1+(6+24g)\log \mu_1 }{48n^{3/2} } + {\rm O}(n^{-2}) \right ),
\EE 
and when $\sigma = \frac{1}{3},$ to
\BE \label{eq:third}
r_n = \mu \left (1 + \frac{ \log \mu_1}{3{n^{2/3}}} + \frac{g}{n} + \frac{\log^2\mu_1}{18n^{4/3}}+ \frac{(2+6g)\log \mu_1 }{18n^{5/3} } + {\rm O}(n^{-2}) \right ).
\EE
and when $\sigma = \frac{1}{4},$ to
\BE \label{eq:fourth}
r_n = \mu \left (1 + \frac{ \log \mu_1}{4{n^{3/4}}} + \frac{g}{n} + \frac{\log^2\mu_1}{32n^{3/2}}+ \frac{(3+8g)\log \mu_1 }{32n^{7/4} } + {\rm O}(n^{-2}) \right ).
\EE

The presence of the term O$(\frac{1}{n^{1-\sigma}})$ in the expression for the ratios above means that a ratio plot against $1/n$ will display curvature, which can be usually be removed by 
 plotting the ratios against $1/n^{1-\sigma},$ with $\sigma =1/2$ or $1/3,$ or $1/4$ etc. A full study of such a situation applied to $Av(1324)$ PAPs can be found in \cite{CGZ18}.
 
 Of course as we saw in the previous section, curvature in the ratio plots can also arise from a pure power-law singularity if the coefficient of the term O$(\frac{1}{n^2})$ in the expression for the ratio is much greater than the coefficient of the term O$(\frac{1}{n}).$ This situation can be identified by plotting the modified ratios, as defined by eqn. (\ref{eqn:modrat}). If these are linear, that is strong evidence of a pure power-law singularity.
 
One can also apply the following test to identify those situations when we have a stretched-exponential singularity. From eqn. (\ref{eq:rn}), we note that, with $s_n = r_n/r_{n-1},$
\BE \label{eqn:div}
n^2(s_n-1) = (\sigma -1)\cdot \sigma \log \mu_1\cdot n^\sigma - g.
\EE(If $\sigma=1/2,$ $g$ must be replaced by $g- \log \mu_1^2/8.$) With a stretched-exponential term, this sequence should diverge with $n.$ In the presence of a pure power-law, the sequence should tend to $-g$ as $n \to \infty.$ We show in Fig.~\ref{fig:A12345div} that in this case the sequence is indeed tending to $-g=7.5$ as expected.

Unfortunately the observation that a ratio plot against $1/n^{1-\sigma}$ will linearise the plot does not provide a sufficiently precise method to estimate the value of $\sigma.$  One can usually distinguish between, say, $\sigma=1/2$ and $\sigma = 1/3$ in this way, but one cannot be much more precise than that. However, as we now show, one can extend the ratio method to provide direct estimates for the value of $\sigma.$

From (\ref{eq:rn}), one sees that 
\BE \label{eq:rsigma}
(r_n/\mu-1) = \sigma \log{\mu_1} \cdot n^{\sigma-1} + O\left ( \frac{1}{n} \right ).
\EE
 Accordingly, a plot of $\log(r_n/\mu-1)$ versus $\log{n}$ should be linear, with gradient $\sigma-1.$ We would expect an estimate of $\sigma$ close to that which linearised the ratio plot.

This log-log plot will usually be visually linear, but the local gradients are changing slowly as $n$ increases. It is therefore worthwhile extrapolating the local gradients. To do this, from (\ref{eq:rsigma}), we form the estimators
\BE \label{eqn:sig1}
{\tilde \sigma}_n = 1+ \frac{\log |r_n/\mu-1|-\log |r_{n-1}/\mu-1|}{\log{n}-\log(n-1)}.
\EE
This can be extrapolated against $1/n^\sigma,$ using any approximate value of $\sigma.$

A second estimator of $\sigma$ follows from eqn. (\ref{eq:an}). Define $$c_n \equiv \log(b_n/\mu^n) \sim \log{C} + \log{\mu_1}\cdot n^\sigma + g\cdot \log{n},$$  then setting
\BE \label{eq:sigma2}
d_n \equiv c_n-c_{n-1} \sim \sigma \log{\mu_1}\cdot n^{\sigma-1} + g/n,
\EE
 a log-log plot of $d_n$ against $n$ should be linear with gradient $\sigma-1.$ Note that if $\sigma$ is closer to zero than to 1, there is likely to be some competition between the two terms in the  expansion.

This way of estimating $\sigma$ requires knowledge of, or at worst a very precise estimate of, the growth constant $\mu.$ While $\mu$ is exactly known in some cases, more generally $\mu$ is not known, and must be estimated, along with all the other critical parameters. In order to estimate $\sigma$ without knowing $\mu,$ we can use one (or both) of the following estimators:

From eqn. (\ref{eq:rn}), it follows that 
\BE \label{eq:sig1}
r_{\sigma_n} \equiv \frac{r_n}{r_{n-1}} \sim 1 + \frac{(\sigma-1)\log{\mu_1}}{n^{2-\sigma}} + {\rm O}(1/n^2),
\EE
so $\sigma$ can be estimated from a plot of $\log(r_{\sigma_n}-1)$ against $\log{n},$ which should have gradient $\sigma-2.$ Again, the local gradients can be calculated and plotted against $1/n^\sigma,$ using any approximate value of $\sigma.$

Another estimator of $\sigma$ when $\mu$ is not known follows from eqn. (\ref{eq:an}),
\BE \label{eq:sig2}
a_{\sigma_n} \equiv \frac{b_n^{1/n}}{b_{n-1}^{1/(n-1)}} \sim 1 + \frac{(\sigma-1)\log{\mu_1}}{n^{2-\sigma}} + {\rm O}(1/n^2),
\EE
so again $\sigma$ can be estimated from a plot of $\log(a_{\sigma_n}-1)$ against $\log{n}.$ Again, estimates of $\sigma$ are found by extrapolating the local gradient  against $1/n^\sigma.$

While these two estimators are equal to leading order, they differ in their higher-order terms. Which of the two is more informative seems to vary from problem to problem. However, we generally use both.

From eqn. (\ref{eq:rn}), if we know (or conjecture) $\mu$ and $\sigma,$ we can use this to estimate $\mu_1,$ as
\BE \label{eqn:mu1}
\left ( \frac{r_n}{\mu} - 1\right )\cdot n^{1-\sigma} \sim \sigma \cdot \log(\mu_1).
\EE

\subsection{Direct fitting.} \label{direct}
Another, perhaps obvious, idea is to try and fit the critical parameters directly to the assumed asymptotic form,  $$b_n \sim C\cdot \mu^n  \cdot n^g$$ in the case of a pure power law singularity, or

$$b_n \sim C\cdot \mu^n \cdot \mu_1^{n^\sigma} \cdot n^g$$ in the case of a stretched exponential singularity.

 For a stretched exponential singularity,
\BE \label{logcan1}
\log {b_n} \sim \log{C} + n \log{\mu} + n^\sigma \log{\mu_1} + g \log {n}.
\EE
So if $\sigma$ is known, or assumed, there are four unknowns in this linear equation. It is then straightforward to solve the linear system
\BE \label{eq:fit1}
\log {b_k} =  c_1k  + c_2 k^\sigma  +c_3 \log {k}+c_4
\EE
 for $k=n-2,\, n-1, \, n, \, n+1$ with $n$ ranging from $3$ to $N-1,$ where $N$ is the highest known power of the series. We refer to this as a {\em 4-point fit.} Then $c_1$ estimates $\log(\mu),$ $c_2$ estimates $\log(\mu_1)$, $c_3$ estimates $g$ and $c_4$ gives estimators of $\log{C}.$ An obvious variation arises in those cases where, say, $\mu$ is known. Then one can solve  
\BE \label{eqn:fit2}
\log(b_k) - k \log{\mu} =  c_1 k^\sigma  +c_2 \log {k} +c_3
\EE
from three successive coefficients, as before increasing the order of the lowest coefficient used by one until one runs out of coefficients. We refer to this as a {\em 3-point fit.}

Alternatively, one can fit the {\em ratios} to 
\BE \label{eqn:ratfit}
r_n = c_1+\frac{c_2}{n^{1-\sigma}} +\frac{c_3}{n}+\frac{c_4}{n^{2-2\sigma}},
\EE
from four successive coefficients, as before increasing the order of the lowest coefficient used  by one until one runs out of coefficients.
Then $c_1$ estimates $\mu,$ $c_2$ estimates $\mu\cdot \sigma \log(\mu_1)$, $c_3$ estimates $\mu\cdot g$ (or $\mu( g-\log^2(\mu_1)/8)$ if $\sigma=1/2$),  and $c_4$ gives estimators of $\mu\cdot \sigma^2 \log^2(\mu_1)/2.$ 

If the value of both $\mu$ and $\sigma$ are known, one can fit three parameters to
\BE \label{eqn:fit3}
r_n-\mu=\frac{c_1}{n^{1-\sigma}}+\frac{c_2}{n} + \frac{c_3}{n^{2-2\sigma}},
\EE
where $c_1$ estimates $\mu \cdot \sigma\log{\mu_1}$  $c_2$ estimates $\mu\cdot g$ and $c_3$ estimates $\mu\cdot \sigma^2 \log^2(\mu_1)/2.$

For the case of a pure power law singularity, we find the idea of direct fitting most useful to estimate the amplitude $C$ when the exponent $g$ is known or conjectured.
In that case one can estimate $\mu$ and $C$ by fitting to
\BE \label{logamp}
\log {b_k} -g\log{n} \sim c_1 n  + c_2 +c_3/n,
\EE
where $c_1$ gives estimators of $\log{\mu}$ and $c_2$ gives estimators of $\log{C}.$

We will apply a number of these techniques in our analysis, below, of the 16 Wilf classes of length-5 PAPS.

\section{Differential approximants}
\label{ana:da}

The generating
functions  of some problems in enumerative combinatorics are sometimes algebraic, such as that for $Av(1342)$ PAPs, sometimes D-finite, such as $Av(12345)$ PAPs,
sometimes differentially algebraic, and sometimes transcendentally transcendental.
The not infrequent occurrence of D-finite solutions was the origin of the method of {\em differential approximants}, a very successful method of series analysis for power-law singularities \cite{G89}.

The basic idea is to approximate a generating function $F(z)$ by solutions
of differential equations with polynomial coefficients. That is to say, by D-finite ODEs. The singular behaviour
of such ODEs is  well documented
(see e.g. \cite{Forsyth02,Ince27}), and the singular points and
exponents are readily calculated from the ODE. 

The key point for series analysis is that even if {\em globally} the function is not describable by a solution
of such a linear ODE (as is frequently the case) one expects that
{\em locally,} in the
vicinity of the (physical) critical points, the generating
function is still well-approximated by a solution of a linear ODE, when the singularity is a generic power law (\ref{generic}).

An $M^{th}$-order differential approximant (DA) to a function $F(z)$  is formed by matching
the coefficients in the polynomials $Q_k(z)$ and $P(z)$ of degree $N_k$ and $L$, respectively,
so that the formal solution of the $M^{th}$-order inhomogeneous ordinary differential equation
\BE \label{eq:ana_DA}
\sum_{k=0}^M Q_{k}(z)(z\frac{{\rm d}}{{\rm d}z})^k \tilde{F}(z) = P(z)
\EE
agrees with the first $N=L+\sum_k (N_k+1)$ series coefficients of $F(z)$. 

Constructing such ODEs only involves
solving systems of linear equations. The function
$\tilde{F}(z)$ thus agrees with the power series expansion of the (generally unknown)
function $F(z)$ up to the first $N$ series expansion coefficients.
We normalise the DA by setting $Q_M(0)=1,$ thus leaving us with $N$ rather
than $N+1$ unknown coefficients to find. The choice of the differential operator $z\frac{{\rm d}}{{\rm d}z}$ in (\ref{eq:ana_DA}) forces the origin to be a regular singular point. The reason for this choice is that most lattice models with holonomic solutions, for example, the free-energy of the two-dimensional Ising model, possess this property. However this is not an essential choice.


From the theory of ODEs, the singularities of $\tilde{F}(z)$ are approximated by zeros
$z_i, \,\, i=1, \ldots , N_M$ of $Q_M(z),$ and the
associated critical exponents $\gamma_i$ are estimated from the indicial equation. If there is only a single root at $z_i$  this is just
\BE \label{eq:ana_indeq1}
\gamma_i=M-1-\frac{Q_{M-1}(z_i)}{z_iQ_M ' (z_i)}.
\EE
Estimates of the critical amplitude $C$ are rather more difficult to make, involving the integration of the differential approximant. For that reason the simple ratio method approach to estimating critical amplitudes is often used, whenever possible taking into account higher-order asymptotic terms \cite{GJ09}.

Details as to which approximants should be used and how the estimates from many approximants are averaged to give a single estimate are given in \cite{GJ09}. Examples of the application of the method can be found in \cite{G15}. In that work, and in this, we reject so-called {\em defective} approximants, typically those that have a spurious singularity closer to the origin than the radius of convergence as estimated from the bulk of the approximants. Another  method sometimes used is to reject outlying approximants, as judged from a histogram of the location of the critical point (i.e. the radius of convergence) given by the DAs. It is usually the case that such distributions are bell-shaped and rather symmetrical, so rejecting approximants beyond two or three standard deviations is a fairly natural thing to do.

\section{Coefficient prediction}
\label{pred}
In \cite{G16} we showed that the ratio method and the method of differential approximants  work serendipitously together in many cases, even when one has stretched exponential behaviour, in which case neither method works particularly well in unmodified form. 

To be more precise, the method of differential approximants (DAs)  produces ODEs which, by construction, have solutions whose series expansions agree term by term with the known coefficients used in their construction. Clearly, such ODEs implicitly define {\em all}  coefficients in the generating function, but if $N$ terms are used in the construction of the ODE, all terms of order $z^{N}$ and beyond will be approximate, unless the exact ODE is discovered, in which case the problem is solved, without recourse to approximate methods.

What we have found is that it is useful to construct a number of DAs that use all available coefficients, and then use these to predict subsequent coefficients. Not surprisingly, if this is done for a large number of approximants, it is found that the predicted coefficients of the term of order $z^n,$ where $n > N,$ agree for the first $k(n)$ digits, where $k$ is a decreasing function of $n.$ We take as the predicted coefficients the mean of those produced by the various DAs, with outliers excluded, and as a measure of accuracy we take the number of digits for which the predicted coefficients agree, or the standard deviation. These two measures of uncertainty are usually in good agreement.

Now it makes no logical sense to use the approximate coefficients as input to the method of differential approximants, as we have used the DAs to obtain these coefficients. However there is no logical objection to using the ({\em approximate}) predicted coefficients as input to the ratio method. Indeed, as the ratio method, in its most primitive form, looks at a graphical plot of the ratios, an accuracy of 1 part in $10^4$ or $10^5$ is sufficient, as errors of this magnitude are graphically unobservable. 

Recall that, in the ratio method one looks at {\em ratios} of successive coefficients. We find that the ratios of the approximate coefficients are predicted with even greater precision than the coefficients themselves by the method of DAs. That is to say, while a particular coefficient and its successor might be predicted with an accuracy of 1 part in $10^p$ for some value of $p$, the {\em ratio} of these successive coefficients is frequently  predicted with significantly greater accuracy (the precision being typically improved by a factor varying between 2 and 20).

The DAs use all the information in the coefficients, and are sensitive to even quite small errors in the coefficients. As an example, in a recent study of some self-avoiding walk series, an error was detected in the twentieth significant digit in a new coefficient, as the DAs were much better converged without the last, new, coefficient. The DAs also require high numerical precision in their calculation. In favourable circumstances, they can give remarkably precise estimates of critical points and critical exponents, by which we mean up to or even beyond 20 significant digits in some cases. Surprisingly perhaps, this can be the case even when the underlying ODE is not D-finite. Of course, the singularity must be of the assumed power-law form.

Ratio methods, and direct fitting methods, by contrast are much more robust. The sort of small error that affects the convergence of DAs would not affect the behaviour of the ratios, or their extrapolants, and would thus be invisible to them. As a consequence, approximate coefficients are just as good as the correct coefficients in such applications, provided they are accurate enough. We re-emphasise that, in the generic situation (\ref{generic}), ratio type methods will rarely give the level of precision in estimating critical parameters that DAs can give. By contrast, the behaviour of ratios can more clearly reveal features of the asymptotics, such as the fact that a singularity is not of power-law type. This is revealed, for example, by curvature of the ratio plots \cite{G15}.

We take, as an example, the OGF for $Av(12453)$ PAPs (see OEIS \cite{OEIS} A116485). This is known to order $x^{38}.$ We will only take the coefficients to order $x^{16}$ and use the method of series extension described above to predict the next 22 ratios, as we can compare them to the exact ratios. The results, based on 3rd order differential approximants, are shown in Table \ref{tab:serpred}. For the first predicted ratio, $r_{18},$ the discrepancy is in the 10th significant digit. For the last predicted ratio, $r_{39}$, the error is in the 5th significant digit. This level of precision is perfectly adequate for ratio analysis.

\begin{table}[htbp]
   \begin{center}
   \begin{tabular}{|l|l|} \hline
   Predicted ratios & Actual ratios \\ \hline
10.463935493 &10.46393544\\
10.654655347& 10.65465504\\ 
10.828226522& 10.82822539\\
10.986854456& 10.98685140\\
11.132386843&11.13238007\\
11.266382111&11.26636895\\
11.390163118&11.39013998\\
11.504857930&11.50482182\\
11.611441483&11.61138359\\
11.710743155&11.71066190\\
11.803496856&11.80338255\\
11.890333733&11.89017822\\
11.971808520&11.97160282\\
12.048402545&12.04814337\\
12.120553112&12.12022972\\
12.188650126& 12.18824275\\
12.252994715&12.25252103\\
12.313939194&12.31336663\\
12.371707700&12.37104982\\
12.426619450&12.42581319\\
12.478784843&12.47787509\\
12.528486946&12.52743256\\
\hline
      \end{tabular}
        \caption{Ratios $r_{18}$ to $r_{39}$ actual and predicted from the coefficients of $Av(12453).$}
   \label{tab:serpred}
   \end{center}
\end{table}

In practice we find that the more exact terms we know, the greater is the number of predicted terms, or ratios that can be predicted. In this study, we typically have 26 or so series terms for each Wilf class. These are usually sufficient to predict about 100 additional ratios to 6-digit accuracy. For the singular class $Av(12453)$ we know 39 terms. In that case we are able to predict 400 further ratios.

In this study, we have extended the ratios of the generating functions of the 15 unknown Wilf classes by typically 100 additional ratios, and have analysed the resulting series by ratio methods.

\section{The permutations}
We have divided the 16 Wilf classes into two sets: those six which we believe have simple power-law asymptotic behaviour, so that $s_n \sim C \cdot \mu^n \cdot n^g,$ and those ten that we believe have stretched-exponential behaviour similar to that exhibited by $Av(1324),$ so that the coefficients  behave asymptotically as $s_n \sim C \cdot \mu^n \cdot \mu_1^{n^\sigma} \cdot n^g,$ where $\sigma =1/2,\,\,\, 1/3$ or $\sigma = 1/4.$ 

We will go through the first examples of each in some detail, to show clearly what is involved in the analysis. Most of the other cases flow similarly {\em mutatis mutandis}, so we give less detail in those cases, except when some distinctive feature warrants further discussion.

Our results are summarised in Tables \ref{tab:growth1} and \ref{tab:growth2} below.

\begin{landscape}
\begin{table}[htbp]
   \begin{center}
   \small
   \begin{tabular}{|l|l|l|l|l|l|} \hline
Pattern&Growth&Lower&Lower&Exponent&Amplitude\\
           &constant $\mu.$&bound, l-c.&bound St.&$g.$&$C.$\\ \hline 
25314&$12.5670 \pm0.0003$&10.8809&12.4622&$-3.214 \pm 0.006$&$0.0164 \pm 0.0001$\\
31524&$12.7344 \pm 0.0003$&11.0042&12.6417&$-2.905 \pm 0.01$&$0.00274 \pm 0.00005$\\
35214&$13.275 \pm0.005$&11.2336&13.1159&$-3.75\pm 0.10$&$0.0145 \pm 0.0015$\\
43251&$13.703 \pm0.001$&11.4821&13.5111&$-4.43 \pm 0.02$&$0.207 \pm 0.005.$\\
34215&$13.945 \pm0.008$&11.6002&13.7131&$-4.67 \pm 0.01$&$0.375 \pm 0.08$\\
12345&$16.0$&14.8735&15.9395&$-7.5 $&275.6\\
\hline
      \end{tabular}
        \caption{Summary of the analysis of the 6 Wilf classes of length-5 PAPs with power-law singularities, ranked by increasing estimates of the growth constant. Lower bound l-c comes from the assumption of log-convexity. Lower bound St. comes from the stronger assumption that the series is a Stieltjes series.}
   \label{tab:growth1}
   \end{center}
\end{table}
\end{landscape}

\begin{landscape}
\begin{table}[htbp]
   \begin{center}
   \small
   \begin{tabular}{|l|l|l|l|l|l|l|l|} \hline
Pattern&Growth&Lower&Lower&Growth&Exponent&Exponent&Amplitude\\
            &constant $\mu.$&bound, l-c.&bound St.&constant $\mu_1.$&$\sigma.$&$g.$&$C.$\\ \hline 

53124&$14.24 \pm0.05$&11.3441&13.5836&--& $0.25 \pm 0.1$&--&--\\
32541&$14.32 \pm 0.04 $&11.5813&13.8447&--&$0.27 \pm 0.07$&--&--\\
35124&$14.54\pm0.03$&11.4025&13.7433&--&$1/4$&--&--\\
31245&$9+4\sqrt{2}=14.6568\ldots$&12.5274&14.3792&$0.27 \pm 0.02$&$1/3$&$-4.7 \pm 0.3$&--\\
42351&$15.10 \pm0.05$&11.5749&14.0314&$0.0012 \pm 0.0006$&$1/3$&$-1.25 \pm 0.25$&$1.5\pm 0.5$\\
42315&$15.40 \pm0.04$&11.8117&14.5633&$0.0001 \pm 0.00005$&$1/4$&--&$90\pm  16$\\
35241&$16.20 \pm0.04$&11.6779&14.6253&$0.00009 \pm 0.00004$&$1/3$&--&$18 \pm 9$\\
53241&$18.66 \pm0.05$&11.9590&15.4445&$0.027 \pm 0.006$&$1/2$&$-4\pm1$&--\\
53421&$19.4092 \pm 0.0003$&12.4079&16.3053&$0.044 \pm 0.002$&$1/2$&$-4.0 \pm 0.5$&$7.4 \pm 2.0$\\
52341&$24.7 \pm 0.3$&12.1992&17.2302&$0.0015$&$1/2$&--&--\\

\hline
      \end{tabular}
        \caption{Summary of the analysis of the 10 Wilf classes of length-5 PAPs with stretched-exponential singularities, ranked by increasing estimates of the growth constant. Lower bound l-c comes from the assumption of log-convexity. Lower bound St. comes from the stronger assumption that the series is a Stieltjes series.}
   \label{tab:growth2}
   \end{center}
\end{table}
\end{landscape}

\subsection{Wilf classes with pure power-law behaviour}

\subsubsection{Av(12345)}
This case is completely solved, \cite{BEGM20}, and the generating function is D-finite. Regev \cite{R80} gave an expression for the asymptotics, $$s_n(12345) \sim C\cdot 16^n/n^{15/2},$$ where $C=3\cdot 2^9/\pi^{3/2}.$  Given the existence of the known solution, it may be thought that there is little point in conducting an analysis. However, it is worthwhile as it shows the strengths and limitations of our methods, admittedly for a special, and possibly benign, case.

\begin{figure}[h!] 
\begin{minipage}[t]{0.45\textwidth} 
\centerline{\includegraphics[width=\textwidth,angle=0]{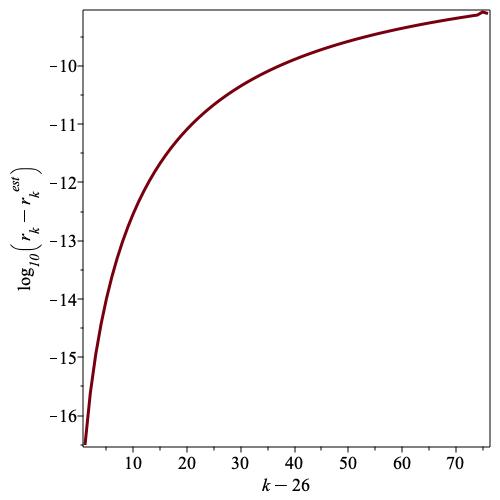} }
\caption{$\log_{10}$ of the difference between the predicted and actual ratios vs. $k$.}
\label{fig:rap12345}
\end{minipage}
\hspace{0.05\textwidth}
\begin{minipage}[t]{0.45\textwidth} 
\centerline{\includegraphics[width=\textwidth,angle=0]{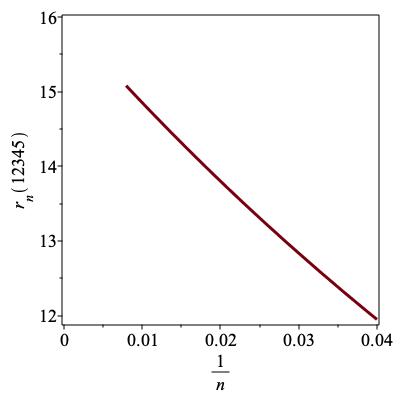}}
\caption{Ratios of $Av(12345)$ vs. $1/n$.}
\label{fig:r112345}
\end{minipage}
\end{figure}

So, to be in accord with our knowledge of most of  the other series, we assume that coefficients are known only to order $x^{25},$  and use these to predict the next 75 ratios. These are increasingly inaccurate with increasing order, but are still useful to the quoted order. In Fig.~\ref{fig:rap12345} we show the base-10 logarithm  of the difference between the predicted ratios $r_k^{\textrm{est}}$ and the true ratios $r_k.$ It can be seen that this difference increases from about $10^{-16.5}$ to $10^{-9}$ as we move from the first predicted ratio to the 75th. (In fact, for several of the other series we use 100 predicted ratios, and in a few cases 200).

Next we carry out an extended ratio analysis. In Fig.~\ref{fig:r112345} we show the predicted ratios plotted against $1/n.$ There is some low order curvature, but the plot is becoming increasingly linear as $n$ increases. The plot is going towards an intercept at $1/n=0$ around 16. 

Following Occam's razor, we assume pure power-law behaviour unless we find evidence to the contrary. For a pure power-law, we expect the linear intercepts $$l_n = n\cdot r_n - (n-1) \cdot r_{n-1}$$ to approach the same limit more closely. As can be readily seen from eqn. (\ref{ratios}), the linear intercepts eliminate the term $O(1/n)$ in the ratios, and in the presence of a pure power-law should vary as $O(1/n^2).$ 

\begin{figure}[h!] 
\begin{minipage}[t]{0.45\textwidth} 
\centerline{\includegraphics[width=\textwidth,angle=0]{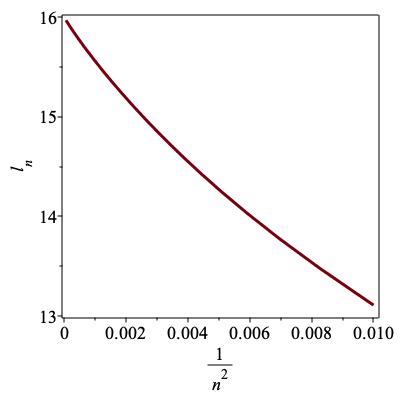} }
\caption{Linear intercepts $l_n$ vs. $1/n^2$.}
\label{fig:r212345}
\end{minipage}
\hspace{0.05\textwidth}
\begin{minipage}[t]{0.45\textwidth} 
\centerline{\includegraphics[width=\textwidth,angle=0]{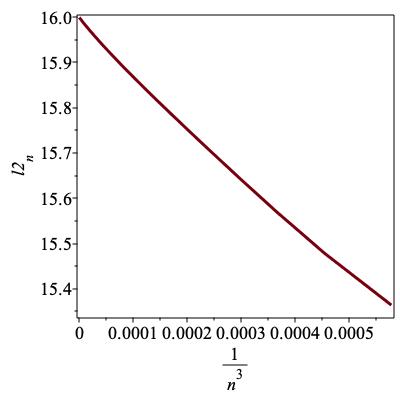}}
\caption{Quadratic intercepts $l2_n$ vs. $1/n^{3}$.}
\label{fig:r312345}
\end{minipage}
\end{figure}

\begin{figure}[h!] 
\begin{minipage}[t]{0.45\textwidth} 
\centerline{\includegraphics[width=\textwidth,angle=0]{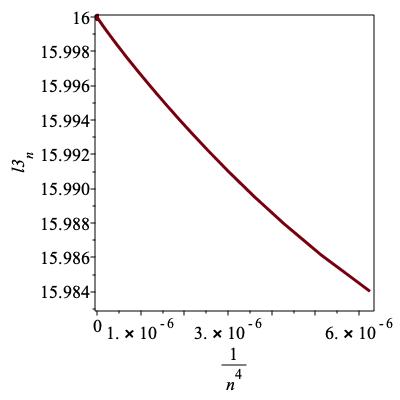} }
\caption{Cubic intercepts $l3_n$ vs. $1/n^4$.}
\label{fig:r412345}
\end{minipage}
\hspace{0.05\textwidth}
\begin{minipage}[t]{0.45\textwidth} 
\centerline{\includegraphics[width=\textwidth,angle=0]{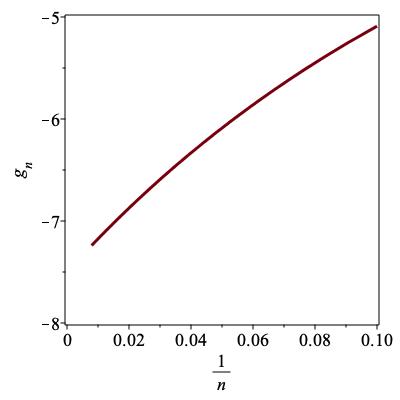} }
\caption{Estimate of exponent $g$, (assuming $\mu$)  vs. $1/n$.}
\label{fig:g112345}
\end{minipage}
\end{figure}

This process can be iterated, and the terms of order $O(1/n^2)$  eliminated by forming the quadratic intercepts $$l2_n = \frac{n^2 \cdot l_n - (n-1)^2 \cdot l_{n-1}}{2n-1},$$ and the cubic intercepts $$l3_n = \frac{n^3 \cdot l2_n - (n-1)^3 \cdot l2_{n-1}}{3n^2-3n+1}.$$
The linear and quadratic and cubic intercepts are shown in  figs. \ref{fig:r212345},  \ref{fig:r312345} and \ref{fig:r412345}  respectively. It can be seen that these are approaching 16.0 more and more precisely\footnote{There is a small amount of ``jitter" in the last few estimates of $l3_n,$ reflecting the error in the (approximate) ratio estimates which gets amplified as one takes higher and higher order differences.}. We feel confident in estimating $\mu = 16.000 \pm 0.005$ on the basis of this study.

\begin{figure}[h!] 
\begin{minipage}[t]{0.45\textwidth} 
\centerline{\includegraphics[width=\textwidth,angle=0]{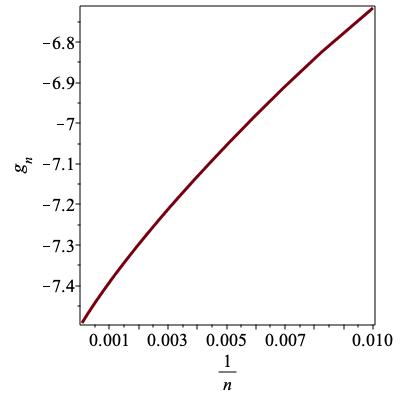} }
\caption{linear extrapolants of exponent $g$, (assuming $\mu$)  vs. $1/n^2$.}
\label{fig:g312345}
\end{minipage}
\hspace{0.05\textwidth}
\begin{minipage}[t]{0.45\textwidth} 
\centerline{\includegraphics[width=\textwidth,angle=0]{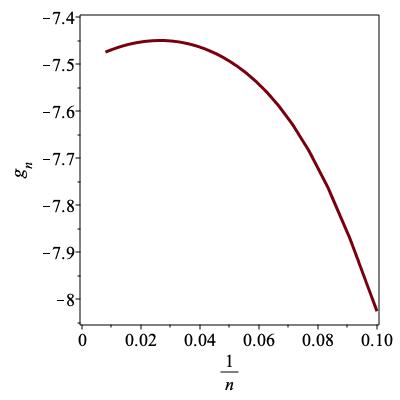}}
\caption{Estimate of exponent $g$ vs. $1/n$ without assuming $\mu.$}
\label{fig:g212345}
\end{minipage}
\end{figure}

Using this estimate of $\mu,$ we can estimate the exponent of the sub-dominant term $n^g.$ From the expression for the ratios (\ref{ratios}), we have that estimators of the exponent $g$ are given by $$g_n = n (r_n/\mu-1) +O(1/n).$$ These are shown in Fig.~\ref{fig:g112345}, and are plausibly approaching a limit $g \approx -7.5.$  As we did with the ratios, these too can be linearly extrapolated by forming the estimates $$g2_n = n\cdot g_n - (n-1) \cdot g_{n-1}.$$ These are shown in Fig.~\ref{fig:g312345}, and are seen to be very convincingly approaching $g = -7.5.$

If $\mu$ is not known, one can still estimate $g$ from the estimators $$g_n =
n^2 (1-r_n/r_{n-1}) +O(1/n).$$ These estimators are shown in Fig.
\ref{fig:g212345}, and are also plausibly approaching a limit of -7.5. The maximum around $n=35$ illustrates the great value of having additional, approximate terms. Without them, one may conclude that the plot was heading toward -7.4 or -7.3.

Of course, one can have power-law behaviour without the ratios behaving as $r_n = \sum_{k=0}^\infty \alpha_k/n^k$ as tacitly assumed in the above analysis. For example, one might have 
\BE \label{eqn:delta}
r_n \sim \mu\left ( 1 + g/n + h/n^{1+\Delta} + j/n^2 + \cdots \right ),
\EE
 where $0 < \Delta < 1.$
 
\begin{figure}[h!] 
 \begin{minipage}[t]{0.45\textwidth} 
\centerline{\includegraphics[width=\textwidth,angle=0]{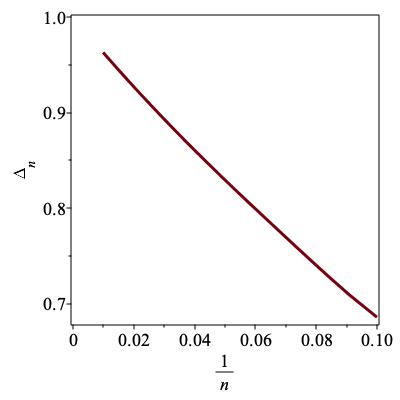} }
\caption{Estimate of exponent $\Delta$ vs. $1/n$.}
\label{fig:Delta12345}
\end{minipage}
\hspace{0.05\textwidth}
\begin{minipage}[t]{0.45\textwidth} 
\centerline{\includegraphics[width=\textwidth,angle=0]{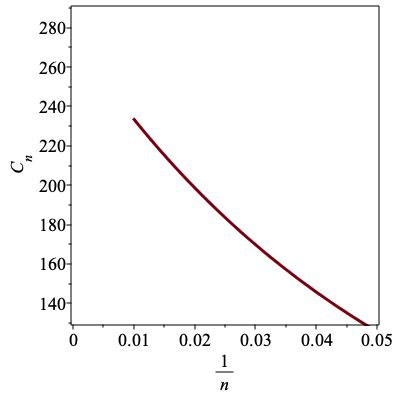}}
\caption{Estimate of amplitude $C$ vs. $1/n$.}
\label{fig:amp12345}
\end{minipage}
\end{figure}

  In the case $Av(12345)$ considered here, we know that we have pure power-law behaviour. Assuming that we do not know, one can try and estimate the value of the exponent $\Delta.$
If one has pure power-law behaviour, one should find $\Delta=1.$ This of course is a necessary, but not sufficient, condition. One might have the more unusual, but by no means impossible behaviour of the ratios $$r_n \sim \mu\left ( 1 + g/n + h/n^2 + j/n^{2+\Delta} +  \cdots\right ).$$ We can only investigate the simpler case with any hope of success. We do this by observing from eqn. (\ref{eqn:delta}) that
$$ \left ( \frac{r_n}{\mu}-1 \right )n - g \sim \frac{h}{n^{\Delta}},$$ so that a log-log plot of the l.h.s. against $n$ should have gradient $ -\Delta.$ In Fig.~\ref{fig:Delta12345} we show the local gradient of the log-log plot, plotted against $1/n,$ which is plausibly going to a limit $\Delta=1.$ This is consistent with the known pure power-law behaviour.

Another way of analysing the series is to fit the expression for the ratios to the available data and solve the resulting system of linear equations, as discussed in Sec. 7.2. Here we will fit the ratios to the assumed form
\BE \label{4fit}
r_n =\mu \left (1 + \frac{g}{n} + \frac{h}{n^2} + \frac{j}{n^3}\right ),
\EE
 using successive quadruplets of ratios $r_{k-2}, r_{k-1}, r_k, r_{k+1},$  with $k$ increasing until we run out of known ratios, to give estimates of the parameters $\mu,$ $g,$ $h,$ and $j.$

\begin{figure}[h!]
\begin{minipage}[t]{0.45\textwidth} 
\centerline{\includegraphics[width=\textwidth,angle=0]{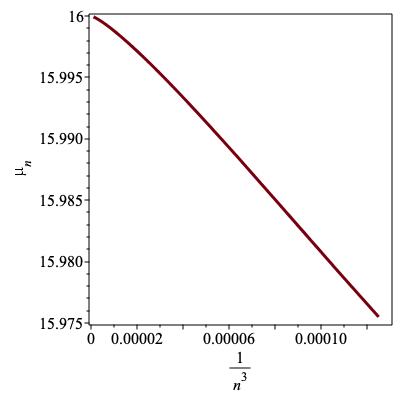} }
\caption{Estimates of $\mu$ vs. $1/n^3$.}
 \label{fig:12345f1}
\end{minipage}
\hspace{0.05\textwidth}
\begin{minipage}[t]{0.45\textwidth} 
\centerline{\includegraphics[width=\textwidth,angle=0]{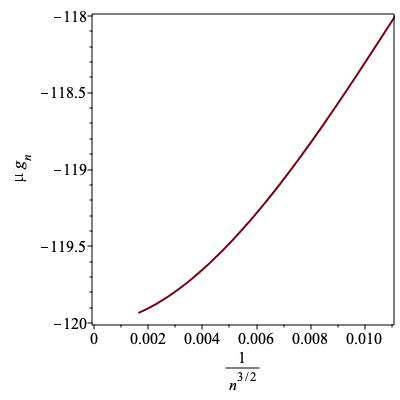}}
\caption{Estimates of $\mu\cdot g$ vs. $1/n^{3/2}$.}
\label{fig:12345f2}
\end{minipage}
\end{figure}

\begin{figure}[h!]
\begin{minipage}[t]{0.45\textwidth}  
\centerline{\includegraphics[width=\textwidth,angle=0]{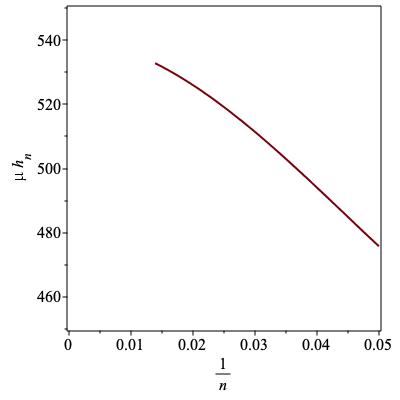} }
\caption{Estimates of $\mu \cdot h$ vs. $1/n$.}
\label{fig:12345f3}
\end{minipage}
\hspace{0.05\textwidth}
\begin{minipage}[t]{0.45\textwidth} 
\centerline{\includegraphics[width=\textwidth,angle=0]{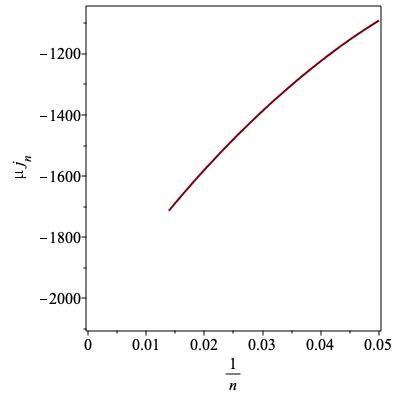}}
\caption{Estimates of $\mu \cdot j$ vs. $1/n$.}
 \label{fig:12345f4}
\end{minipage}
\end{figure}

The results of this fitting are shown in the four figures figs. \ref{fig:12345f1}, \ref{fig:12345f2}, \ref{fig:12345f3}, \ref{fig:12345f4}. From the first figure, one would feel confident estimating $\mu$ at 16.0. exactly. Similarly, from the second figure, the estimate $\mu \cdot g \approx -120$ looks quite compelling, from which we conclude $g \approx -7.5,$ as before. From the next two figures we estimate $\mu \cdot h \approx 540,$ and $\mu \cdot j \approx -2100,$ so that $h \approx 33.75,$ and $j \approx 131.25.$

If we insert these parameters into eqn. (\ref{4fit}), one predicts $r_{100}=14.8519.$ The correct value is $14.851999203\ldots,$ so this asymptotic form, based on only 25 known coefficients, is pleasingly accurate. 

To emphasise what has been achieved here, we took just 26 coefficients, predicted 75 further ratios using differential approximants, fitted these 100 ratios to the assumed asymptotic form, (having already produced good numerical evidence for the overall structure), and derived an asymptotic expression for the ratios which is accurate to 6 significant digits at order 100.


Finally, we can estimate the amplitude $C$ in the expression for the asymptotic form of the coefficients, $s_n \sim C \cdot 16^n/n^{7.5}$ by forming simple estimators of $C$ defined by $$C_n \equiv s_n\cdot n^{7.5}/16^n.$$ These are shown Fig.~\ref{fig:amp12345}, plotted against $1/n$ and are going to a limit around $C \approx 275.$ This numerical estimate could be refined, but is known exactly (see OEIS A047889) as $C=3\cdot 2^9/\pi^{3/2}=275.8458\ldots .$

Alternatively, if we knew, or conjectured, the exponent $g=-7.5,$ we could use direct fitting, as described in eqn. (\ref{logamp}) to estimate $C$ (and indeed $\mu$). We did this and obtained the estimates $\mu \approx 16.0000$ and $C \approx 275.6.$ These are of course in complete agreement with the estimates obtained from the ratios.

If we didn't know the exact result, we would conclude  that $\mu = 16.000 \pm 0.005, \,\, g = -7.5 \pm 0.2,$ and $C \approx 275.6.$ That is to say $$s_n(12345)\sim C \mu^n \cdot n^g.$$ Of course, we know that the central estimates of $\mu$ and $g$ are exact, and that $C=3\cdot 2^9/\pi^{3/2}=275.8458\ldots .$ 

\subsubsection{Lower bounds}
It is provably the case \cite{BEGM20} that  all possible Hankel determinants constructed from the coefficients of $Av(12345)$ are positive, and monotonically increasing with the size of the matrix. As the coefficients of $Av(12345)$ therefore form a Stieltjes moment sequence,  log-convexity of the coefficients follows, as discussed above, 
and so the ratios provide an increasing sequence of lower bounds. This gives the bound $\mu(12345) \ge 14.8735$ from the first 100 coefficients.

If one constructs the continued fraction representation from the exact coefficients, (see Theorem 1 above), then the terms $\alpha_0 \ldots \alpha_{100}$  defined in Theorem 1 can be
used to construct stronger bounds $(\sqrt{\alpha_n} + \sqrt{\alpha_{n-1}})^2.$ Using the first 100 exact coefficients, we find the strong lower bound $\mu(25314) \ge 15.9395.$\
  
\begin{figure}[h!]
\centerline{\includegraphics[width=4in,angle=0]{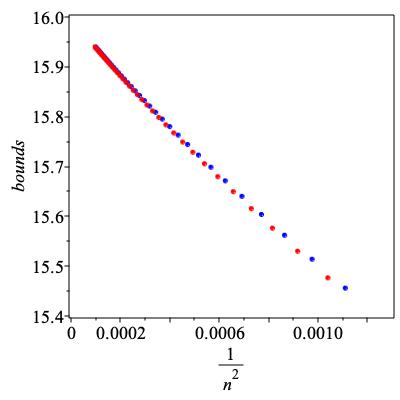} }
\caption{ Stieltjes bounds vs. $1/n^2$ for $Av(12345).$}
 \label{fig:b12345}
 \end{figure}

  The sequence of lower bounds can also be extrapolated against $1/n^2,$ shown
  in Fig.~\ref{fig:b12345}, and results in an estimate of $\mu$ consistent with, but less precise than, that obtained by the direct analysis of the original, extended series.

We  constructed all possible Hankel determinants from all the known (actual) coefficients for the remaining 15 unsolved Wilf classes, and observed that they too are all positive, and monotonically increasing with the size of the matrix, leading us to the conjecture that they too form a Stieltjes moment sequence. We will assume this in the analysis of the remaining 15 Wilf classes, and so provide (conjectured) lower bounds, repeating the above analysis {\em mutatis mutandis.}

The analysis of the remaining five 
Wilf classes with power-law singularities follows along similar lines, and is given in Appendix A.

\subsection{Wilf classes with stretched-exponential behaviour}

As discussed above in Sec. \ref{sec:se}, the (naive) hallmark of such asymptotic behaviour is non-linearity -- more precisely curvature -- in the ratio plots. However, as seen in the discussion of the PAP $Av(12345),$ a power-law singularity can have a ratio plot exhibiting curvature if the magnitude of the coefficient of the term O$(1/n^2)$ is significantly greater than that of the term O$(1/n).$
In that case, we formed quadratic intercepts to eliminate the O$(1/n^2)$ term, and the linearity of the ratio plot, expected for a power-law singularity, became clear.

In the case of stretched exponential singularities, the ratios behave as 
\begin{multline} \label{eq:rna}
r_n = \mu \left (1 + \frac{\sigma \log \mu_1}{n^{1-\sigma}} + \frac{g}{n} + \frac{\sigma^2 \log^2 \mu_1}{2n^{2-2\sigma}} + \frac {(\sigma-\sigma^2)\log \mu_1+2g\sigma \log \mu_1}{2n^{2-\sigma}} \right . \\
 \left . {}+ \frac{\sigma^3 \log^3 \mu_1}{6n^{3-3\sigma}} +{\rm O}(n^{2\sigma-3}) + {\rm O}(n^{-2}) \right ).
\end{multline}

Any ratio plot will manifest competition between the term of order O$(\frac{1}{n^{1-\sigma}} )$ and that of order O$(\frac{1}{n} ).$ We can eliminate the latter term by constructing the linear intercepts 
\begin{multline} \label{eq:lna}
l_n \equiv n \cdot r_n - (n-1) \cdot r_{n-1} \sim \mu \left (1 + \frac{\sigma(\sigma-1)\log \mu_1}{n^{1-\sigma}}  + {\rm O} \left (n^{2\sigma-2} \right ) + {\rm O} \left ( n^{\sigma-2} \right ) \right . \\
 \left . {}+ {\rm O} \left ( {n^{3\sigma-3}} \right )+{\rm O}(n^{2\sigma-3}) + {\rm O}(n^{-2}) \right ).
\end{multline}

Then one can eliminate the term of order O$(1/n^2)$ by forming the quadratic intercepts

\begin{multline} \label{eq:l2na}
l2_n \equiv \frac{ n^2 \cdot l_n - (n-1)^2 \cdot l_{n-1}}{2n-1} \sim \mu \left (1 + \frac{\sigma^2(\sigma-1)\log \mu_1}{2n^{1-\sigma}}  + {\rm O} \left (n^{2\sigma-2} \right )  \right . \\
 \left . {}+ {\rm O} \left ( n^{\sigma-2} \right ) + {\rm O} \left ( {n^{3\sigma-3}} \right )+{\rm O}(n^{2\sigma-3})  \right ).
\end{multline}

So these quadratic intercepts eliminate both the O$(1/n)$  and the O$(1/n^2)$ terms. Plotting $l2_n$ against $1/n^{1-\sigma}$ should give a linear plot at the correct value of $\sigma.$ We find this to
be a more reliable indicator of the (approximate) value of the exponent $\sigma$ than the simpler criterion of linearity of the ratios.

\subsection{ Av(12453)}

This series is known up to, and including, terms of order $x^{38},$ given in the OEIS as sequence A116485. It is exceptional in two respects. Firstly in that so many coefficients are known, due to Biers-Ariel, who showed  \cite{B19} that this pattern has properties that allowed him to write a particularly efficient algorithm for its enumeration. 

Secondly, the growth constant is known exactly \cite{B05}, as $\mu=9+4\sqrt{2} = 14.65685\ldots$.

\begin{figure}[h!] 
\begin{minipage}[t]{0.45\textwidth} 
\centerline{\includegraphics[width=\textwidth,angle=0]{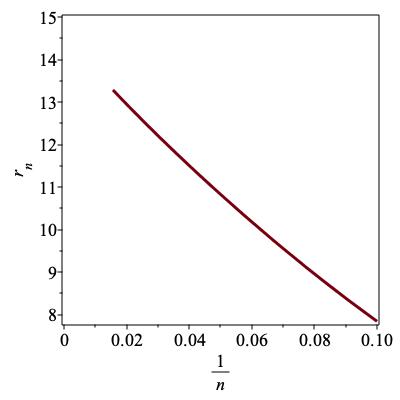} }
\caption{Ratios of $Av(12453)$ vs. $1/n$.}
\label{fig:r112453}
\end{minipage}
\hspace{0.05\textwidth}
\begin{minipage}[t]{0.45\textwidth} 
\centerline{\includegraphics[width=\textwidth,angle=0]{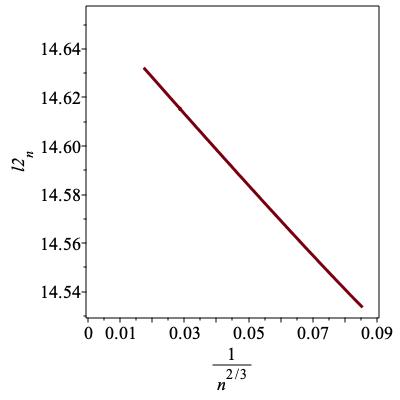}}
\caption{Quadratic intercepts of ratios vs. $1/n^{2/3}$.}
\label{fig:l212453}
\end{minipage}
\end{figure}

With this longer series we have been able to extend the ratio sequence, and coefficient sequence by 400 further terms. For maximum precision, we only use the first 200 terms in the subsequent analysis. 
The ratio plot, showing the ratios plotted against $1/n$ is shown in Fig.~\ref{fig:r112453}, and exhibits convex (w.r.t. the $x$-axis) curvature, which is not disappearing as $n$ increases, unlike the situation arising for pure power-law singularities. Plotting the modified ratios $r2_n,$ defined in eqn (\ref{eqn:modrat}), also fails to linearise the ratio plot. This is strong evidence for a stretched-exponential singularity.

However, when the ratios are plotted against $n^{-\theta},$ where $\theta$ is in the range $(0.75--0.85)$ the plot is visually linear, and appears to extrapolate to the known value of $\mu.$ However, when we plot the quadratic intercepts, as defined in eqn. (\ref{eq:l2na}), which eliminate the effects of the  O$(1/n)$ and  O$(1/n^2)$ terms in the ratios, the plot is actually linearised when plotted against $n^{-\theta}$ with $\theta \approx 2/3,$ as shown in Fig.~\ref{fig:l212453}.

 To better estimate $\sigma,$ we consider two ways of doing so, assuming $\mu$ is known, as it is. From 
 eqn.(\ref{eq:rsigma}), one sees that a plot of $\log(r_n/\mu-1)$ versus $\log{n}$ should be linear, with gradient $\sigma-1.$ So we estimate $\sigma$ from the gradient of this log-log plot. We would expect an estimate of $\sigma$ close to that which linearised the quadratic intercepts plot. 
As terms of order O$(1/n)$  affect the linearity of the ratio plots, they will also have a similar effect here. So to eliminate that effect, we instead plot $\log(l_n/\mu-1)$ versus $\log{n}.$
 
 A second estimator of $\sigma$ follows from eqn.(\ref{eq:sigma2}),
so that a plot of $\log(d_n)$ against $\log{n}$ should be linear, again with gradient $\sigma-1.$

\begin{figure}[h!] 
\begin{minipage}[t]{0.45\textwidth} 
\centerline{\includegraphics[width=\textwidth,angle=0]{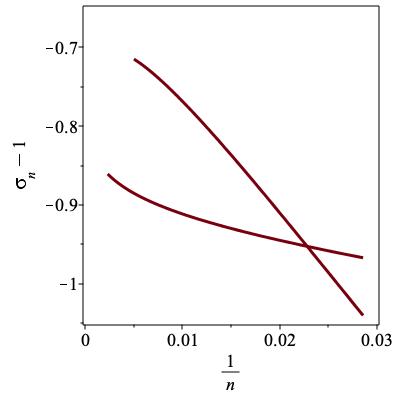}}
\caption{Estimate of exponent $\sigma-1$ vs. $1/n$. First method, upper curve (as $1/n \to 0$), second method, lower curve.}
\label{fig:sig12453}
\end{minipage}
\hspace{0.05\textwidth}
\begin{minipage}[t]{0.45\textwidth} 
\centerline{\includegraphics[width=\textwidth,angle=0]{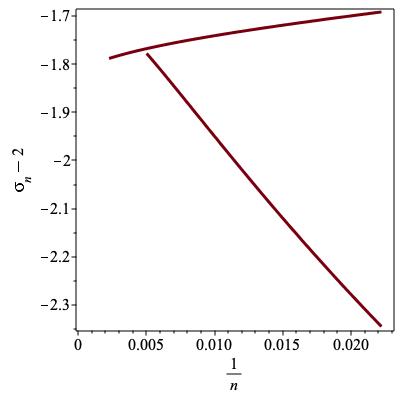} }
\caption{Estimate of exponent $\sigma-2$ vs. $1/n$ without knowing $\mu.$ First method, upper curve, second method, lower curve.}
\label{fig:sig212453}
\end{minipage}
\end{figure}

Both estimators will usually provide visually linear log-log plots, but the local gradients are changing as $n$ increases. One must extrapolate the local gradients. To do this we form the estimators (in the first case)
\BE
{\tilde \sigma}_n -1= \frac{\log |l_n/\mu-1|-\log |l_{n-1}/\mu-1|}{\log{n}-\log(n-1)}.
\EE
This can be extrapolated against $1/n.$ 
The result of doing this is shown in Fig.~\ref{fig:sig12453}, upper curve, as $1/n$ approaches zero. 
The local gradient of the second estimator, given by eqn. (\ref{eq:sigma2}), can be similarly extrapolated, and similarly plotted against $1/n.$ The result of doing so is shown in Fig.~\ref{fig:sig12453}, lower curve.

It can be seen that the second estimator displays considerable curvature, and is difficult to extrapolate, beyond saying that $\sigma > 0.15,$ but the first estimator, which we expect to be more believable as it eliminates the effect of the O$(1/n)$ term and the O$(1/n^2)$ term in the ratios, is plausibly going to $\sigma \approx 1/3.$

An estimator of $\sigma$ formed without knowledge of the growth constant $\mu$ follows from eqn.(\ref{eq:sig1}), so that a plot of $\log(r_{\sigma_n}-1)$ against $\log{n},$  should have gradient $\sigma-2.$  However, for the same reasons as above, a more reliable estimator should follow by replacing $r_n$ by $l_n$ in eqn. (\ref{eq:sig1}), and this results in the top curve in \ref{fig:sig212453}.

Another estimator of $\sigma$ when $\mu$ is not known is given in eqn. (\ref{eq:sig2}),
so  $\sigma-2$ can be estimated from a plot of $\log(a_{\sigma_n}-1)$ against $\log{n}.$ 

Again, estimates of $\sigma$ are found by extrapolating the local gradient  against $1/n^\sigma.$
The result of doing this is shown in Fig.~\ref{fig:sig212453}, (lower curve). 

These estimates are, unsurprisingly, less precise than the estimates formed knowing $\mu.$ Indeed, the second estimator is seemingly going to a value inconsistent with the other estimators, which are consistent with our previous estimate $\sigma \approx 1/3.$ We ascribe this to the strong effect of the $O(1/n)$ and O$(1/n^2)$ term in the ratios for this pattern.

\begin{figure}[h!] 
\begin{minipage}[t]{0.45\textwidth} 
\centerline{\includegraphics[width=\textwidth,angle=0]{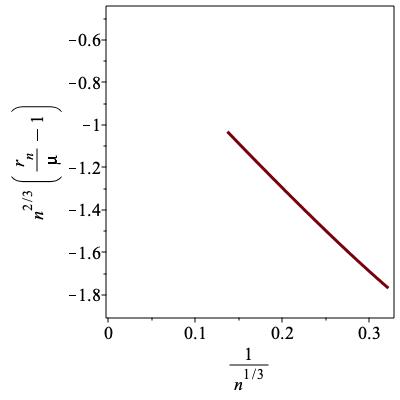}}
\caption{Estimate of  $\mu_1$ vs. $1/n^{1/3}$.}
\label{fig:mu112453}
\end{minipage}
\hspace{0.05\textwidth}
\begin{minipage}[t]{0.45\textwidth} 
\centerline{\includegraphics[width=\textwidth,angle=0]{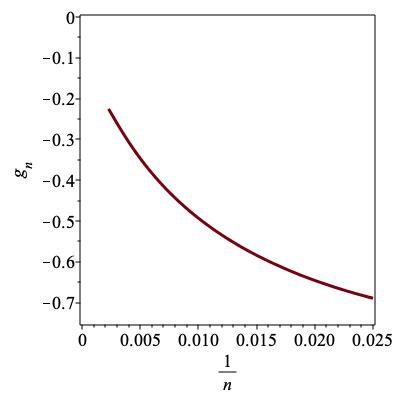} }
\caption{Estimators of exponent $g$ vs. $1/n$ for $Av(12453).$}
\label{fig:g12453}
\end{minipage}
\end{figure}

If we accept $\sigma = 1/3$ as the most likely exact value, we can estimate the sub-dominant growth constant $\mu_1.$ From eqn (\ref{eq:rsigma}) it follows that $$ \left ( \frac{r_n}{\mu} -1 \right ) n^{1-\sigma} \sim \sigma \cdot \log{\mu_1}+\frac{g}{n^\sigma},$$ and we show in Fig.~\ref{fig:mu112453} a plot of the l.h.s. against $1/n^{1/3}.$ From this curve, we estimate $\sigma \cdot \log{\mu_1} \approx -0.45,$ so $\mu_1 \approx 0.26.$ 

\begin{figure}[h!] 
\begin{minipage}[t]{0.45\textwidth} 
\centerline{\includegraphics[width=\textwidth,angle=0]{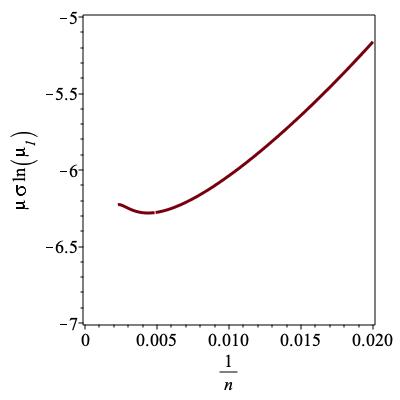}}
\caption{Estimate of  $\mu\cdot \sigma \log \mu_1$ vs. $1/n$.}
\label{fig:param1}
\end{minipage}
\hspace{0.05\textwidth}
\begin{minipage}[t]{0.45\textwidth} 
\centerline{\includegraphics[width=\textwidth,angle=0]{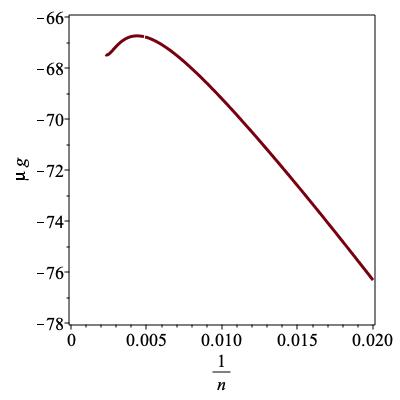} }
\caption{Estimators of exponent $\mu \cdot g$ vs. $1/n.$ }
\label{fig:param2}
\end{minipage}
\end{figure}

An alternative, and more precise estimate of the value of $\mu_1$ can be made as follows: Knowing the value of $\mu,$ and assuming $\sigma=1/3,$ one can fit the remaining parameters, to the expression for the ratios, as described in eqn. (\ref{eqn:fit3}). We show the results of this in figs. \ref{fig:param1} and \ref{fig:param2}.

The first parameter, which we estimate to be $-6.2 \pm 0.2$ gives $\mu\cdot \sigma \log{\mu_1}.$ Therefore $\mu_1 = 0.28 \pm 0.01,$ in quite good agreement with the less-precise estimate $0.26$ given above.  

The second parameter, which we estimate to be $-68.5\pm 1.5$ gives $\mu\cdot g.$ Therefore $g=-4.67 \pm 0.1.$  The third parameter, not shown, is harder to estimate, but gives us no further information, as the only unknown involved is $\log \mu_1,$ which we have
 already estimated from the first parameter.

Alternatively, one can fit to the expression for the logarithm of the coefficients, as given in eqn. (\ref{eqn:fit2}), and doing this we find $\log \mu_1 \approx -1.25$ and $g \approx -4.7.$ So $\mu_1 \approx 0.27.$ With this degree of variation in the estimates  of $\mu_1$ and $g,$ it is not possible to reliably estimate the amplitude $C.$

 We therefore conclude that $s_n(12453)\sim C \cdot \mu^n \cdot \mu_1^{n^{\sigma}}  \cdot n^g,$ with $\mu = 9+4\sqrt{2} = 14.65685\ldots, \,\, \sigma \approx 1/3, \,\, \mu_1 = 0.27 \pm 0.02,$ and $g = -4.7 \pm 0.3.$ We give no estimate of the amplitude $C.$ 
 
 These parameter values explain some of the difficulties we encountered in this analysis. The ratios behave as $$r_n \sim \mu \left ( 1 -\frac{0.44}{n^{2/3}} -\frac{4.7}{n} + \cdots \right ).$$ One sees that the coefficient of the O$(n^{-2/3})$ term is about 1/10 of that of the O$(1/n)$ term, which explains why we had to eliminate the O$(1/n)$ and O$(1/n^2)$ terms to be able to focus on those terms that arose from the stretched-exponential.
 That said, we still find the presence of a stretched-exponential term for this pattern mildly surprising, as several other PAPs which clearly exhibit stretched-exponential behaviour  (42315, 35241, 53241, 53421, 52341) contain a sub-pattern that qualitatively looks like $1324,$ while this pattern looks like $3124,$ which has a power-law singularity.

 \begin{figure}[h!] 
\centerline{\includegraphics[width=3in,angle=0]{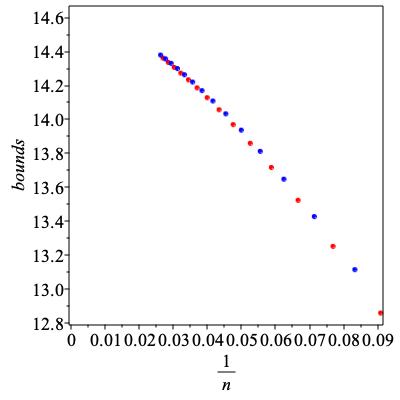}}
\caption{  Stieltjes bounds vs. $1/n$ for $Av(12453).$.}
\label{fig:b12453}
\end{figure}

\subsubsection{Lower bounds}
While there is little point in obtaining bounds in this case where the growth constant is exactly known, it is perhaps of interest to see how close the bounds are to the exact value.
Assuming that the coefficients form a Stieltjes moment sequence, the consequent log-convexity gives the bound $\mu(12453) \ge 12.5274.$ If the 400 predicted coefficients are believed, this improves the bound to $\mu(12453) \ge 14.4008.$ From the continued-fraction representation, we obtain the lower bound $14.3792.$ 
In Fig.~\ref{fig:b12453} we show these Stieltjes bounds plotted against the appropriate power of $n.$ It can be seen they are extrapolating 
 to the top-leftmost corner of the plot, which corresponds to $\mu.$


The analysis of the remaining nine
Wilf classes with stretched-exponential singularities follows along similar lines, and is given in Appendix B.

\section{Conclusion}
This is the first of a series of papers we propose, analysing the sixteen length-5 classical pattern-avoiding permutations. The analysis given here is based on series expansions. We are also developing new Monte Carlo algorithms to provide further information, such as the shape of a typical permutation, as well as providing an alternative route to investigating the growth constants.

Here we have extended the number of known coefficients for fourteen of the sixteen classes, see table \ref{tab:terms1} and \ref{tab:terms2}  . Using sequence extension and a variety of methods of series analysis  we have estimated the growth constant of all classes, and  have estimate the sub-dominant power-law term associated with the exponential growth, or the stretched-exponential exponent. Repeating the {\em caveat} made in the introduction, the only singularity types we are considering here
are pure power-law and stretched exponentials, motivated of course by the fact that these are the only singularity types we have encountered for PAPs of shorter length.
However, if there is another singularity type, or even stretched exponentials with additional logarithms terms, we are not testing for that.

In six of the sixteen classes cases we found the familiar power-law behaviour, so that the coefficients behave like $s_n \sim C \mu^n n^g,$ while in the remaining
ten  cases we find a stretched exponential as the sub-dominant term, so that the coefficients behave like $s_n \sim C \mu^n \mu_1^{n^\sigma} n^g,$ where $0 < \sigma < 1.$ 
 
 We have also classified the 120 possible permutations into the 16 distinct classes, see table \ref{tab:classes}.

We give lower bounds to the growth constant in all cases, based on the belief,  and in one case a proof, that all 16 Wilf-class generating function coefficients can be represented as Stieltjes moment sequences.

The numerical estimates of growth constants and sub-dominant terms are given in tables \ref{tab:growth1} and \ref{tab:growth2}.

\section{Acknowledgements}
AJG would like to thank the ARC Centre of Excellence for Mathematical and Statistical Frontiers (ACEMS) for support. This research was undertaken using the Research Computing Services facilities hosted at the University of Melbourne.

\section {Appendix A}
In this appendix we give details of the analysis of the remaining five Wilf classes with a power-law singularity.

\subsection{ Av(25314)}
The series was known up to, and including, terms of order $x^{16},$ given in the OEIS as sequence A256195. We have extended the known series by ten further terms, and the sequence of approximate coefficients and ratios by 100 further terms,
using the method of series extension described above, based on third order differential approximants. 

\begin{figure}[h!]
\begin{minipage}[t]{0.45\textwidth} 
\centerline{\includegraphics[width=\textwidth,angle=0]{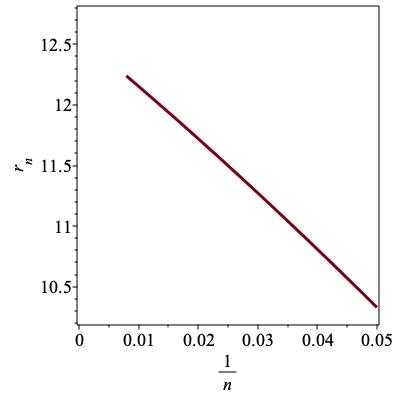} }
\caption{Ratios of $Av(25314)$ vs. $1/n$.}
 \label{fig:r25314}
\end{minipage}
\hspace{0.05\textwidth}
\begin{minipage}[t]{0.45\textwidth} 
\centerline{\includegraphics[width=\textwidth,angle=0]{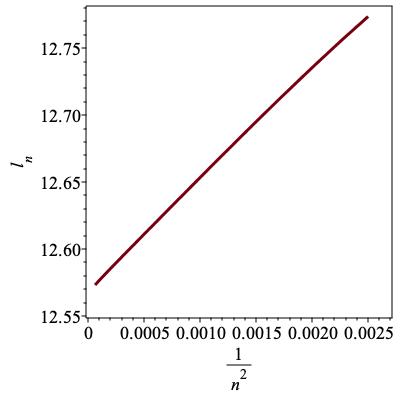}}
\caption{Linear intercepts of $Av(25314)$ vs. $1/n^2$.}
\label{fig:l25314}
\end{minipage}
\end{figure}

\begin{figure}[h!]
\begin{minipage}[t]{0.45\textwidth}  
\centerline{\includegraphics[width=\textwidth,angle=0]{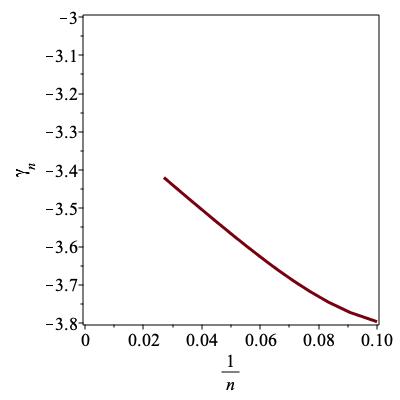} }
\caption{Exponent estimates of $Av(25314)$ assuming $\mu=12.5670$.}
\label{fig:g25314}
\end{minipage}
\hspace{0.05\textwidth}
\begin{minipage}[t]{0.45\textwidth} 
\centerline{\includegraphics[width=\textwidth,angle=0]{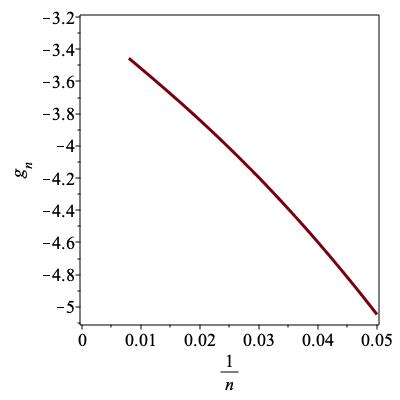}}
\caption{Exponent estimate of $Av(25314),$ independent of $\mu.$}
 \label{fig:d25314}
\end{minipage}
\end{figure}

A plot of the ratios against $1/n$ is shown in Fig.~\ref{fig:r25314},
and we estimate the extrapolated limit at $1/n = 0$ to be $\mu = 12.6 \pm 0.1.$ In Fig.~\ref{fig:l25314} we show the linear intercepts, from which we give the more precise estimate
$\mu=12.57 \pm 0.01.$ With quadratic intercepts we can sharpen this slightly to $12.5670 \pm 0.0003.$ Using the central estimate $\mu=12.5670,$ we show estimates of the exponent $g$  in Fig.~\ref{fig:g25314} which we can linearly extrapolate, and give the estimate $g = -3.213 \pm 0.005.$

To increase confidence that we have the same pure power-law behaviour as prevails for $Av(12345),$ we estimated the exponent $\Delta,$ defined in eqn. (\ref{eqn:delta}), and found it to be extremely close to $1,$ consistent with pure power-law behaviour. 

In Fig.~\ref{fig:d25314}, we show a plot of estimators of the exponent $g$ which is independent of any estimate of $\mu.$ This extrapolates to a value consistent with the estimate
just made assuming a value for $\mu.$  Linearly extrapolating this plot also gives the estimate $g = -3.213 \pm 0.005.$


\begin{figure}[h!]
\begin{minipage}[t]{0.45\textwidth} 
\centerline{\includegraphics[width=\textwidth,angle=0]{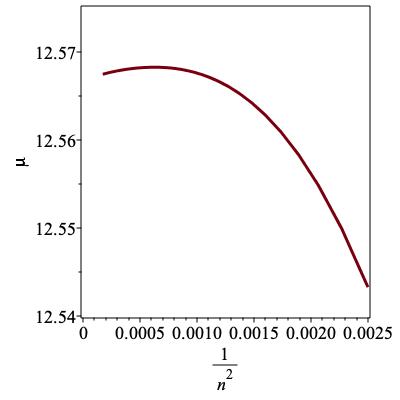} }
\caption{Estimate of growth constant $\mu$ vs. $1/n^2$.}
 \label{fig:mu25314}
\end{minipage}
\hspace{0.05\textwidth}
\begin{minipage}[t]{0.45\textwidth} 
\centerline{\includegraphics[width=\textwidth,angle=0]{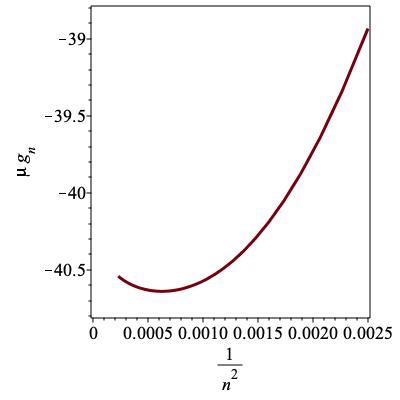}}
\caption{Estimate of $\mu\cdot g$ vs. $1/n^2$.}
\label{fig:mug25314}
\end{minipage}
\end{figure}

We also fitted the ratio sequence to the assumed asymptotic form eqn. (\ref{4fit}), and show the results in figs. \ref{fig:mu25314} and \ref{fig:mug25314} for $\mu$ and $\mu\cdot g$ respectively.
Our previous estimate $\mu \approx 12.5670$ is well-supported, while $\mu \cdot g \approx -40.4,$ so $g \approx -3.215,$  both being consistent with our previous estimates.  

We expect that the exponent is rational. The closest simple rational number is $3 \frac{3}{14} = 3.21425\cdots,$ and we offer this as our best guess as to the exact value. Using the central estimate of $\mu$ and the conjectured exponent, we can estimate the amplitude $C,$ as we did for $Av(12345)$ in the previous subsection. In this way we find $C = 0.0164 \pm 0.0001.$

As we did for the analysis of  $Av(12345),$ using the conjectured value of the exponent $g=-3 \frac{3}{14}$ in this case, we can use direct fitting, as described in eqn. (\ref{logamp}) to estimate $C$ (and indeed $\mu$). We did this and obtained the estimates $\mu \approx 12.5670$ and $C \approx 0.0164.$ These are of course in complete agreement with the estimates obtained from the ratios.

 We therefore conclude $\mu = 12.5670 \pm 0.0003, \,\, g = -3.214 \pm 0.006,$  and $C  = 0.0164 \pm 0.0001.$ That is to say $$s_n(25314)\sim C \mu^n \cdot n^g,$$ where our best guess is  that $g=-3 \frac{3}{14}$ exactly.  \\
 
\subsubsection{Lower bounds}
Assuming that the Hankel determinants remain positive, then the ratios of the coefficients provide an increasing sequence of lower bounds. This gives the bound $\mu(25314) \ge 10.8809 .$ If the 100 predicted ratios are accepted, this improves the bound to $\mu(25314) \ge 12.240.$

If one constructs the continued fraction representation from the exact coefficients, (see Theorem 1 above), then the terms $\alpha_0 \ldots \alpha_{26}$  defined
in Theorem 1 can be
used to construct stronger bounds $(\sqrt{\alpha_n} + \sqrt{\alpha_{n-1}})^2.$ In this way we obtain the bound $\mu(25314) \ge 12.46223.$

\begin{figure}[h!]
\begin{minipage}[t]{0.45\textwidth} 
\centerline{\includegraphics[width=\textwidth,angle=0]{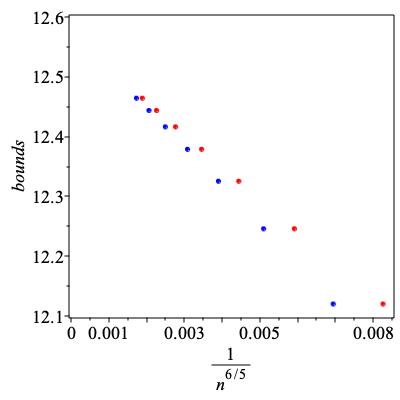} }
\caption{ Stieltjes bounds vs. $1/n^2$ for $Av(25314).$}
 \label{fig:b25314}
 \end{minipage}
\hspace{0.05\textwidth}
\begin{minipage}[t]{0.45\textwidth} 
\centerline{\includegraphics[width=\textwidth,angle=0]{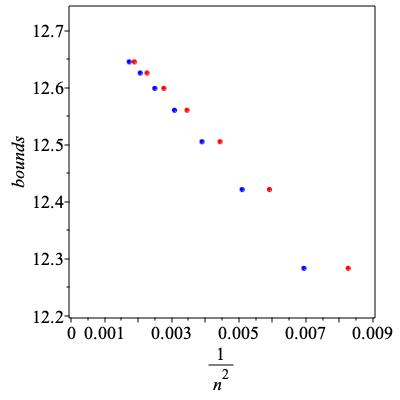}}
\caption{ Stieltjes bounds vs. $1/n^2$ for $Av(31524).$}
\label{fig:b31524}
\end{minipage}
\end{figure}

The sequence of lower bounds can also be extrapolated against $1/n^2,$ shown in Fig.~\ref{fig:b25314}, and results in an estimate of $\mu$ consistent with, but less precise than, that obtained by the direct analysis of the original, extended series.

\subsection{ Av(31524)}
This series is now known up to, and including, terms of order $x^{24},$ given in the OEIS as sequence A256196 to order $x^{16}.$ We have extended the sequence of ratios and coefficients by 100 further terms, as described above. The analysis parallels that just described above for the case of $Av(25314).$  \\

From the ratios, we estimate $\mu=12.75 \pm 0.1.$ From the linear intercepts we give the more precise estimate
$\mu=12.74 \pm 0.01.$ Quadratic intercepts sharpen this to $\mu=12.734 \pm 0.001,$ and cubic intercepts give $\mu=12.7344 \pm 0.0003.$ Using the central estimate $\mu=12.7344,$ we estimate $g = -2.905 \pm 0.01.$
Our estimate of the exponent $g$ independent of any estimate of $\mu$  extrapolates to a value consistent with this estimate. As we consider it likely that the exponents are simple rationals, our best guess is $g=-2\frac{9}{10}.$
We also made a 4-point fit to the ratios, and a 5-point fit, as defined in Sec. \ref{direct}. These confirmed our estimates, above,  of $\mu$ and $g,$ and provided abundant evidence of a pure power-law singularity. 
 As in the previous pattern, we estimated the amplitude $C,$ and found $C = 0.00274 \pm 0.00005.$

As we did for the analysis of  $Av(12345),$ using the conjectured value of the exponent $g=-2\frac{9}{10}$ in this case, we can use direct fitting, as described in eqn. (\ref{logamp}) to estimate $C$ (and indeed $\mu$). We did this and obtained the estimates $\mu \approx 12.7344$ and $C \approx 0.000275.$ These are of course in complete agreement with the estimates obtained from the ratios.

 We therefore conclude $\mu = 12.7344 \pm 0.0003, \,\, g = -2.905 \pm 0.01$ and $C = 0.00274 \pm 0.00005.$ That is to say $$s_n(31524)\sim C \mu^n \cdot n^g,$$ where our best guess is that $g=-2 \frac{9}{10}$ exactly.  \\

\subsubsection{Lower bounds}
Assuming the observed positivity of the Hankel matrices persists,
it follows that the coefficients form a Stieltjes moment sequence. The consequent log-convexity gives the bound $\mu(31524) \ge 10.9278.$ If the 100 predicted coefficients are accepted, this improves the bound to $\mu(31524) \ge 12.4169.$

Constructing the continued fraction representation from the exact coefficients, (see Theorem 1 above), as in the previous pattern, then one finds  $\mu(31524) \ge 12.6417$

 The sequence of lower bounds can also be extrapolated against $1/n^2,$ shown in Fig.~\ref{fig:b31524}, and results in an estimate
of $\mu$ consistent with, but less precise than, that obtained by the direct analysis of the original, extended series.

\subsection{ Av(35214)}
This series was known up to, and including, terms of order $x^{14},$ given in the OEIS as sequence A256197. We have extended the series by twelve terms. We have also extended the sequence of ratios and coefficients by 100 further ratios and coefficients, as described above, though the series is not quite as well-behaved, so our confidence bounds are wider.
The analysis parallels that just described above for the case of $Av(31524).$  \\

From the ratios, we estimate $\mu=13.3 \pm 0.1.$ From the linear and quadratic intercepts we give the more precise estimate
$\mu=13.28 \pm 0.01$ and $13.275 \pm 0.005$ respectively. Using the central estimate $\mu=13.275,$ we estimate $g = -3.75 \pm 0.10.$
Our estimate of the exponent $g$ independent of any estimate of $\mu$  extrapolates to a value consistent with this estimate.
We consider it likely that the exponents are simple rationals, so our best guess is  $g=-3.75.$
We also made  4-point and 5-point fits to the ratios, as above. These confirmed our estimates, above,  of $\mu$ and $g,$ and provided abundant evidence of a pure power-law singularity. As for the previous pattern, we estimated the amplitude $C,$ and found $C = 0.0145 \pm 0.0015.$ Note that this estimate assumes the central values for the parameters $\mu$ and $g.$

As we did for the analysis of  $Av(12345),$ using the conjectured value of the exponent $g=-3.75$ in this case, we can use direct fitting, as described in eqn. (\ref{logamp}) to estimate $C$ (and indeed $\mu$). We did this and obtained the estimates $\mu \approx 13.275$ and $C \approx 0.0143.$ These are of course in complete agreement, within quoted uncertainties, with the estimates obtained from the ratios.

 We therefore conclude $\mu = 13.275 \pm 0.005, \,\, g = -3.75 \pm 0.10,$ and $C=0.0145 \pm 0.0015.$ That is to say $$s_n(35214)\sim C \mu^n \cdot n^g,$$ where our best guess is that $g=-3.75$ exactly.  \\
\\

\subsubsection{Lower bounds}
Assuming that the coefficients form a Stieltjes moment sequence, the consequent log-convexity gives the bound $\mu(35214) \ge 11.2336.$ If the 100 predicted coefficients are accepted, this improves the bound to $\mu(31524) \ge 12.8561.$

Constructing the continued fraction representation from the exact coefficients, (see Theorem 1 above), as above, one finds  $\mu(35214) \ge 13.1159.$ 

The sequence of lower bounds can also be extrapolated against $1/n^2,$ and results in an estimate
of $\mu$ consistent with that obtained by the direct analysis of the original, extended series.

\subsection{ Av(43251)}

This series was known up to, and including, terms of order $x^{14},$ given in the OEIS as sequence A256203. We have extended the series by twelve further terms. We have also extended the sequence of ratios by 100 further ratios, as described above.
The analysis parallels those cases already discussed. \\

From the ratios, we estimate $\mu=13.70 \pm 0.02.$ From the linear and quadratic intercepts we give the more precise estimate
$\mu=13.703 \pm 0.001.$   Using the central estimate $\mu=13.703,$ we estimate $g = -4.43 \pm 0.02.$
Our estimate of the exponent $g$ independent of any estimate of $\mu$  extrapolates to a value consistent with this estimate.
We consider it likely that the exponents are simple rationals, and in this case appears to be  $g=-4 \frac{3}{7}.$
We also made  4-point and 5-point fit to the ratios, as above. These confirmed our estimates, above,  of both $\mu$ and $g,$ and provided abundant evidence of a pure power-law singularity. As for the previous pattern, we estimated the amplitude $C,$ and found $C = 0.208 \pm 0.005.$\\

As we did for the analysis of  $Av(12345),$ using the conjectured value of the exponent $g=-4 \frac{3}{7}$ in this case, we can use direct fitting, as described in eqn. (\ref{logamp}) to estimate $C$ (and indeed $\mu$). We did this and obtained the estimates $\mu \approx 13.703$ and $C = 0.206 \pm 0.005$ These are of course in complete agreement, within quoted uncertainties, with the estimates obtained from the ratios.

 We therefore conclude $\mu = 13.703 \pm 0.001, \,\, g = -4.43 \pm 0.02,$ and $C =  0.207 \pm 0.005.$ That is to say $$s_n(43251)\sim C \mu^n \cdot n^g,$$ where our best guess is that $g=-4 \frac{3}{7}$ exactly.  \\


\subsubsection{Lower bounds}
Assuming that the coefficients form a Stieltjes moment sequence, the consequent log-convexity gives the bound $\mu(43251) \ge 11.4821.$ If the 200 predicted coefficients are believed, this improves the bound to $\mu(43251) \ge 13.4367.$

Constructing the continued fraction representation from the exact coefficients, (see Theorem 1 above), as above, then one finds  $\mu(43251) \ge 13.5111.$

 This sequence of lower bounds can also be extrapolated against $1/n^2,$ and results in an estimate
of $\mu$ consistent with that obtained by the direct analysis of the original, extended series.

\subsection{ Av(34215)}

This series was known up to, and including, terms of order $x^{14},$ given in the OEIS as sequence A256205. We have extended the series by twelve further terms. We have also extended the sequence of ratios by 200 further ratios, as described above. We were able to extend the series by 200 ratios rather than the extension by 100 ratios of the preceding sequences simply because the computed error in the coefficients increased more slowly with increasing order of the coefficients in this case than in other cases. We do not know why.
The analysis parallels that described above, though the series is slightly less well-behaved, so the uncertainties in our estimates are greater.  \\

From the ratios, we estimate $\mu=13.97 \pm 0.05.$ From the linear intercepts we give the more precise estimate
$\mu=13.95 \pm 0.01,$ and from quadratic fits and 4-point fits we make a slightly lower estimate, $\mu=13.945 \pm 0.008.$ 
Using the central estimate $\mu=13.945,$ we estimate $g = -4.67 \pm 0.01.$  Our estimate of the exponent $g$ independent of any estimate of $\mu$  extrapolates to a similar value, though less precisely.

Assuming $g=-4\frac{2}{3}$ exactly, then allowing for the error estimate in $\mu,$ we find $C = 0.375 \pm 0.08,$ where the bulk of the uncertainty stems from the uncertainty in $\mu.$

 We therefore conclude $\mu = 13.945 \pm 0.008, \,\, g = -4.67 \pm 0.01,\,\, C=0.375 \pm 0.08.$ That is to say $$s_n(34215)\sim C \mu^n \cdot n^g.$$    Our best guess for the exact value of the exponent $g$ is $g=-4\frac{2}{3}.$ Note that the growth rate is smaller than that of $Av(53124),$ even though all known coefficients of $Av(34215)$ are greater than those of $Av(53124).$ For sufficiently high order, this pattern must of course reverse.

\subsubsection{Lower bounds}
Assuming that the coefficients form a Stieltjes moment sequence, the consequent log-convexity gives the bound $\mu(34215) \ge11.6002.$ If the 200 predicted coefficients are believed, this improves the bound to $\mu(34215) \ge 13.6611.$

Constructing the continued fraction representation from the exact coefficients, (see Theorem 1 above), as above, one finds  $\mu(34215) \ge 13.7131.$\\ This sequence of lower bounds can also be extrapolated against $1/n^2,$ and results in an estimate
of $\mu$ consistent with, but less precise than, that obtained by the direct analysis of the original, extended series.

\section {Appendix B}

In this appendix we give details of the analysis of the remaining nine Wilf classes with a stretched-exponential singularity.

\subsection{ Av(53124)}

This series was known up to, and including, terms of order $x^{14},$ given in the OEIS as sequence A256199. We have extended the series by eleven further terms. 

The analysis parallels that described above for the case of $Av(12453),$ though this series does not behave as nicely as does the series $Av(12453),$ presumably because we have significantly fewer terms. Accordingly, we have been able to extend the number of approximate coefficients, and ratios, by 100, rather than 400 in the previous case.

\begin{figure}[h!]
\begin{minipage}[t]{0.45\textwidth} 
\centerline{\includegraphics[width=\textwidth,angle=0]{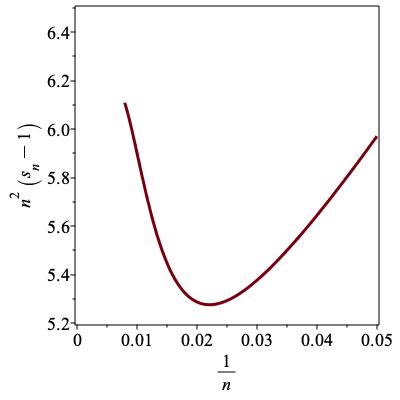} }
\caption{$n^2(s_n-1)$ vs. $1/n,$ showing divergence.}
 \label{fig:d53124}
\end{minipage}
\hspace{0.05\textwidth}
\begin{minipage}[t]{0.45\textwidth} 
\centerline{\includegraphics[width=\textwidth,angle=0]{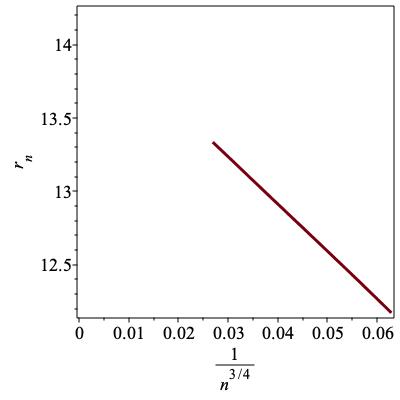}}
\caption{Ratios of $Av(53124)$ vs. $1/n^{3/4}.$}
\label{fig:r253124}
\end{minipage}
\end{figure}

 The ratios, when plotted against $1/n$ display some curvature which is not disappearing as $n$ increases. This allows us to estimate $\mu > 14.0.$ Curvature in this ratio plot is usually a hallmark of a stretched-exponential singularity, as discussed above. We investigate this further by checking if the sequence $n^2(r_n/r_{n-1}-1)$ diverges with $n,$ as explained in the discussion around eqn. (\ref{eqn:div}). The relevant plot is shown in Fig.~\ref{fig:d53124}, and is clearly diverging as $n$ increases. 
 
 If we plot the ratios against $1/n^{3/4},$ see Fig.~\ref{fig:r253124} the plot is visually linear. This would imply a stretched-exponential term of the form $\mu_1^{n^{1/4}},$ which is not an exponent we have previously encountered. Extrapolating this gives the estimate $14.2 < \mu < 14.3.$

The permutation pattern 53124 is a simple decreasing sequence  (531) followed by an increasing sequence, (24), so there is nothing there that would suggest a stretched exponential term. The previous pattern we have studied $Av(12453)$ is similar, in that it is an increasing sequence (1245) followed by a decreasing sequence (53). Most of the sequences we study here that display clear evidence of a stretched-exponential term contain the pattern 1324, whereas this does not.

\begin{figure}[h!]
\begin{minipage}[t]{0.45\textwidth} 
\centerline{\includegraphics[width=\textwidth,angle=0]{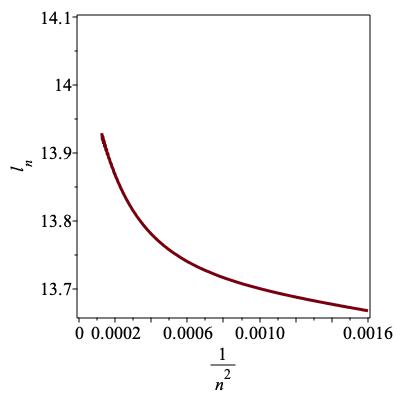}}
\caption{Linear intercepts of $Av(53124)$ vs. $1/n^2$.}
\label{fig:l53124}
\end{minipage}
\hspace{0.05\textwidth}
\begin{minipage}[t]{0.45\textwidth} 
\centerline{\includegraphics[width=\textwidth,angle=0]{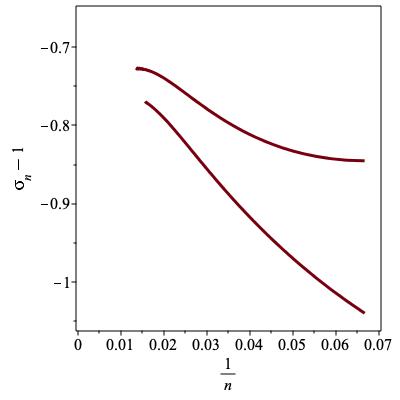} }
\caption{Estimate of exponent $\sigma,$ assuming $\mu$ vs. $1/n.$}
 \label{fig:sig153124}
 \end{minipage}
\end{figure}

If we had a pure power law, linear intercepts would eliminate the $O(1/n)$ term in the expression for the ratios, and the linear intercepts would behave as $l_n \sim \mu(1+c/n^2).$ A plot of the linear intercepts against $1/n^2$ is shown in Fig.~\ref{fig:l53124}, and it is clear that this is far from linear! This is strong evidence against a pure power-law singularity. However plotting these linear intercepts against $1/n^{3/4},$ as appropriate, allows us to estimate $\mu = 14.24 \pm 0.05.$

\begin{figure}[h!]
\begin{minipage}[t]{0.45\textwidth} 
\centerline{\includegraphics[width=\textwidth,angle=0]{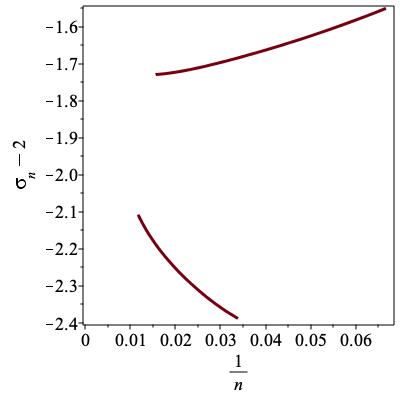}}
\caption{Estimate of exponent $\sigma,$ independent of $\mu$ vs. $1/n.$}
\label{fig:sig253124}
\end{minipage}
\hspace{0.05\textwidth}
\begin{minipage}[t]{0.45\textwidth} 
\centerline{\includegraphics[width=\textwidth,angle=0]{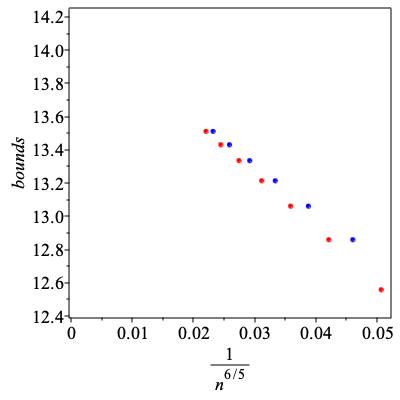} }
\caption{ Stieltjes bounds vs. $1/n^{6/5}$.}
 \label{fig:b53124}
\end{minipage}
\end{figure}

As in the previous pattern, we can form estimates for the exponent $\sigma$ assuming the value of $\mu,$ and, less precisely, without assuming $\mu.$ The result of doing this is shown in figs \ref{fig:sig153124} and \ref{fig:sig253124} respectively, where the upper and lower curves in both figures represent the results of using the two different methods for estimating $\sigma,$ as explained in the previous subsection. It can be seen that these plots lend support to the estimate $\sigma=1/4,$ originally conjectured based on linearity of the ratio plots.

We tried a variety of other methods to estimate the various parameters, as discussed in the previous subsection, but these were inconclusive. Accordingly, while $\sigma=1/4$ is our candidate for the most likely exact value, we cannot totally rule out $\sigma=1/3,$ as found for $Av(12453).$  This is reflected in the quoted error bar, which should be interpreted as a confidence limit.

 The uncertainty in our estimate of $\mu$ prevents us from estimating the stretched-exponential growth constant $\mu_1$ or the exponent $g$ with useful precision. We therefore conclude $\mu = 14.24 \pm 0.05, \,\, \sigma=0.25 \pm 0.1.$ That is to say $$s_n(53124)\sim C \mu^n \cdot \mu_1^{n^\sigma}\cdot n^g,$$ with $\mu_1$ $C$ and $g$ unknown.\\


\subsubsection{Lower bounds}
Assuming that the coefficients form a Stieltjes moment sequence, the consequent log-convexity gives the bound $\mu(53124) \ge 11.3441.$ If the 100 predicted coefficients are accepted, this improves the bound to $\mu(53124) \ge 13.3481.$

Constructing the continued fraction representation from the exact coefficients, (see Theorem 1 above), as above,  one finds  $\mu(53124) \ge 13.5836.$

In Fig.~\ref{fig:b53124} we show these Stieltjes bounds plotted against the appropriate power of $n.$ It can be seen that they can be extrapolated approximately linearly to the top-leftmost point of the plot, at the estimated value $\mu=14.24.$

\subsection{ Av(32541)}

This series was also known up to, and including, terms of order $x^{14},$ given in the OEIS as sequence A256204. We have extended the series by thirteen further terms. We have also extended the sequence of ratios by 100 further ratios, as described above. 
The analysis parallels that just described above for the case of $Av(35124).$  \\

\begin{figure}[h!]
\begin{minipage}[t]{0.45\textwidth} 
\centerline{\includegraphics[width=\textwidth,angle=0]{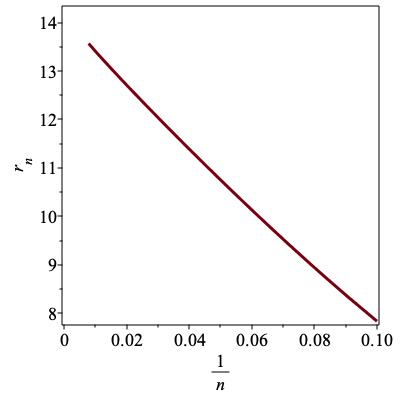}}
\caption{Ratios of $Av(32541)$ vs. $1/n$.}
\label{fig:r132541}
\end{minipage}
\hspace{0.05\textwidth}
\begin{minipage}[t]{0.45\textwidth} 
\centerline{\includegraphics[width=\textwidth,angle=0]{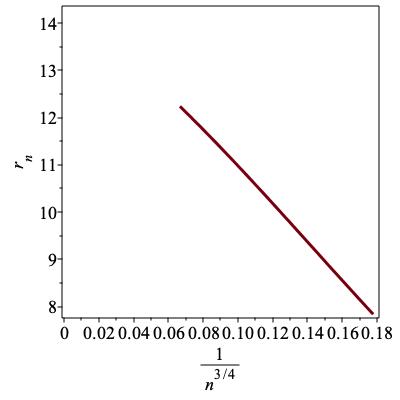} }
\caption{ Ratios of $Av(53124)$ vs. $1/n^{3/4}$..}
 \label{fig:r232541}
 \end{minipage}
\end{figure}

In Fig.~\ref{fig:r132541} we show the ratios plotted against $1/n.$ This plot is slightly concave, and can be crudely extrapolated to give the estimate $\mu \approx 14.3.$ If we plot the ratios against $1/n^{3/4},$ see Fig.~\ref{fig:r232541} the plot is visually linear. This would imply a stretched-exponential term of the form $\mu_1^{n^{1/4}}.$  Extrapolating this gives the estimate $\mu \approx 14.35.$

As in the previous pattern, we can form estimates for the exponent $\sigma$ assuming the value of $\mu,$ and, less precisely, without assuming $\mu.$ The result of doing this is shown in figs \ref{fig:sig135124} and \ref{fig:sig235124} respectively, where the upper and lower curves in both figures represent the results of using two different methods for estimating $\sigma,$ as explained in the previous subsection. It can be seen that these plots are somewhat inconclusive, suggesting estimates of $\sigma$ in the range $(1/5,1/3),$ but with $\sigma=1/4$ as the central estimate. Thus we take this as our most likely exact value, but quote our estimate as $\sigma= 0.27 \pm 0.07.$

Assuming $\sigma=1/4,$ we can fit to the expression for the logarithm of coefficients as described around eqn. (\ref{logcan1}) to estimate the other parameters. In this way we estimate $\mu \approx 14.30,$ $\log \mu_1 \approx -5,$ and $g \approx -2.$ 

With $\sigma=1/3,$ we find $\mu \approx 14.35,$  $\log \mu_1 \approx -2.2,$ and $g \approx -3.$ 

We have also estimated these parameters by fitting the {\em ratios} to the expected form, as explained in the description of eqn. (\ref{eqn:ratfit}). This produces similar estimates.

\begin{figure}[h!]
\begin{minipage}[t]{0.45\textwidth} 
\centerline{\includegraphics[width=\textwidth,angle=0]{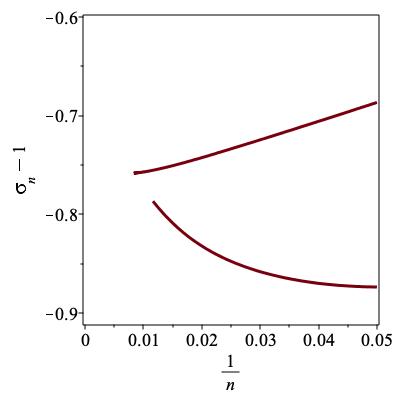} }
\caption{Estimate of exponent $\sigma,$ assuming $\mu$ vs. $1/n.$}
 \label{fig:sig132541}
\end{minipage}
\hspace{0.05\textwidth}
\begin{minipage}[t]{0.45\textwidth} 
\centerline{\includegraphics[width=\textwidth,angle=0]{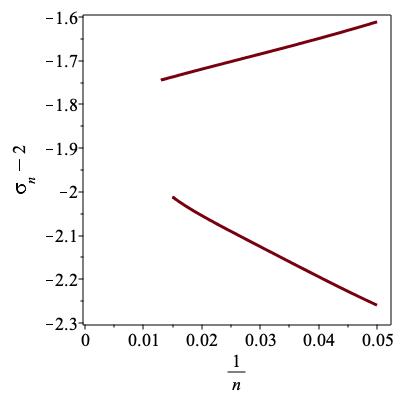}}
\caption{Estimate of exponent $\sigma,$ independent of $\mu$ vs. $1/n.$}
\label{fig:sig232541}
\end{minipage}
\end{figure}

As in the previous pattern considered, the uncertainty in the estimate of $\mu$ prevents us from estimating the stretched-exponential growth constant $\mu_1$ or the exponent $g$ with useful precision, as their estimation depends sensitively on the value of the exponent $\sigma.$ 

We therefore conclude $\mu = 14.32\pm 0.04, \,\, \sigma = 0.27 \pm 0.07,$ with $\sigma=1/4$ or $1/3$ as the most likely values. That is to say $$s_n(32541)\sim C \mu^n \cdot \mu_1^{n^\sigma} \cdot n^g.$$ If  $\sigma = 1/4,$ then we estimate $\mu_1 \approx 0.007,$ and $g \approx -2.$  If  $\sigma = 1/3,$ then we estimate $\mu_1 \approx 0.11,$ and $g \approx -3.$

\subsubsection{Lower bounds}
Assuming that the coefficients form a Stieltjes moment sequence, the consequent log-convexity gives the bound $\mu(32541) \ge 11.5813.$ If the 100 predicted coefficients are believed, this improves the bound to $\mu(32541) \ge 13.5669.$

Constructing the continued fraction representation from the exact coefficients, (see Theorem 1 above), as above, one finds  $\mu(32541) \ge 13.8447.$

 This sequence of lower bounds can also be extrapolated against $1/n^{6/5},$ shown in Fig.~\ref{fig:b32541} and results in an estimate
of $\mu$ consistent with that obtained by the direct analysis of the original, extended series.

\subsection{ Av(35124)}
This series was known up to, and including, terms of order $x^{14},$ given in the OEIS as sequence A256198. We have extended the series by thirteen further terms. We have also extended the sequence of ratios by 100 further ratios, as described above. The analysis parallels that just described above for the case of $Av(32541).$  \\

\begin{figure}[h!]
\begin{minipage}[t]{0.45\textwidth} 
\centerline{\includegraphics[width=\textwidth,angle=0]{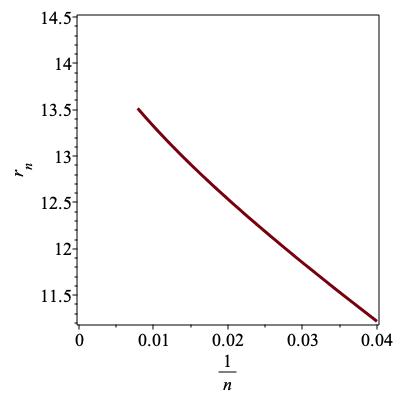}}
\caption{Ratios of $Av(35124)$ vs. $1/n$..}
\label{fig:r35124}
\end{minipage}
\hspace{0.05\textwidth}
\begin{minipage}[t]{0.45\textwidth} 
\centerline{\includegraphics[width=\textwidth,angle=0]{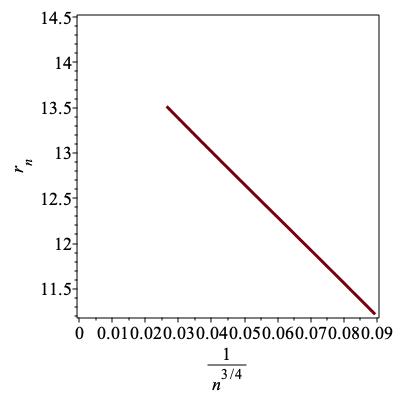} }
\caption{ Ratios of $Av(35124)$ vs. $1/n^{3/4}$.}
 \label{fig:r235124}
 \end{minipage}
\end{figure}

In Fig.~\ref{fig:r35124} we show the ratios plotted against $1/n,$ which appears to be concave. If we plot the ratios against $1/n^{3/4},$ see Fig.~\ref{fig:r235124} the plot is visually linear. This would imply a stretched-exponential term of the form $\mu_1^{n^{1/4}},$ just as observed for $Av(53124).$ Extrapolating this gives the estimate $\mu \approx 14.50.$

If we had a pure power law, linear intercepts would eliminate the $O(1/n)$ term in the expression for the ratios, and the linear intercepts would behave as $l_n \sim \mu(1+c/n^2).$ A plot of the linear intercepts against $1/n^2$ displays considerable curvature, which is further evidence against simple power-law behaviour. However plotting the linear intercepts against $1/n^{3/4},$ as appropriate, allows us to estimate $\mu = 14.54
 \pm 0.03.$

\begin{figure}[h!]
\begin{minipage}[t]{0.45\textwidth} 
\centerline{\includegraphics[width=\textwidth,angle=0]{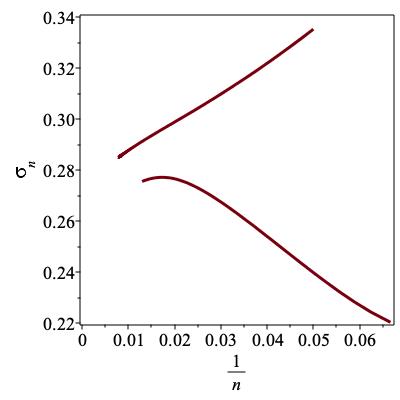} }
\caption{Estimate of exponent $\sigma,$ assuming $\mu$ vs. $1/n.$}
 \label{fig:sig135124}
\end{minipage}
\hspace{0.05\textwidth}
\begin{minipage}[t]{0.45\textwidth} 
\centerline{\includegraphics[width=\textwidth,angle=0]{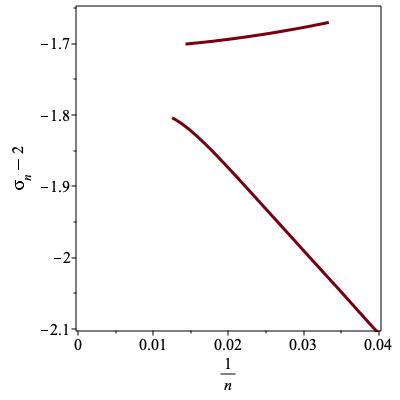}}
\caption{Estimate of exponent $\sigma,$ independent of $\mu$ vs. $1/n.$}
\label{fig:sig235124}
\end{minipage}
\end{figure}

As in the previous pattern, we can form estimates for the exponent $\sigma$ assuming the value of $\mu,$ and, less precisely, without assuming $\mu.$ The result of doing this is shown in figs \ref{fig:sig135124} and \ref{fig:sig235124} respectively, where the upper and lower curves in both figures represent the results of using the two different methods for estimating $\sigma,$ as explained in the previous subsection. It can be seen that these plots lend support to the estimate $\sigma=1/4,$ based on linearity of the ratio plots.

As in the previous pattern considered, the uncertainty in the estimate of $\mu$ prevents us from estimating the stretched-exponential growth constant $\mu_1$ or the exponent $g$ with useful precision. We therefore conclude $\mu = 14.54\pm 0.03, \,\, \sigma=\frac{1}{4}.$ That is to say $$s_n(35124)\sim C \mu^n \cdot \mu_1^{n^\sigma} \cdot n^g.$$ Unlike the previous two patterns, the value $\sigma=1/4$ seems much less equivocal than in those cases, so we do not quote error estimates in this case, but remark that while this is our most favoured value, it is not impossible that $\sigma$ is as large as $1/3.$\\

 \subsubsection{Lower bounds}
Assuming that the coefficients form a Stieltjes moment sequence, the consequent log-convexity gives the bound $\mu(35124) \ge 11.4025.$ If the 100 predicted coefficients are accepted, this improves the bound to $\mu(35124) \ge 13.5150.$

Constructing the continued fraction representation from the exact coefficients, (see Theorem 1 above), as above, then  one finds  $\mu(35124) \ge 13.7433.$

\begin{figure}[h!]
\begin{minipage}[t]{0.45\textwidth} 
\centerline{\includegraphics[width=\textwidth,angle=0]{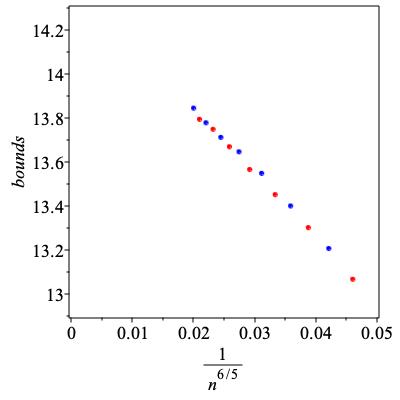}}
\caption{Stieltjes bounds for $Av(32541)$ vs. $1/n^{6/5}$}
\label{fig:b32541}
\end{minipage}
\hspace{0.05\textwidth}
\begin{minipage}[t]{0.45\textwidth} 
\centerline{\includegraphics[width=\textwidth,angle=0]{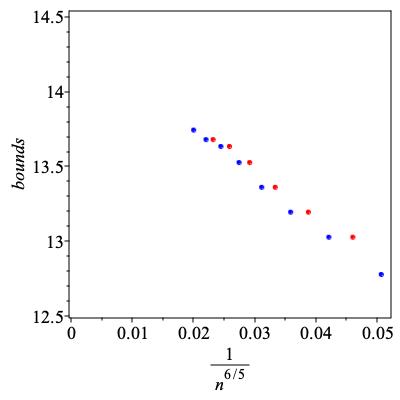} }
\caption{Stieltjes bounds for $Av(35124)$ vs. $1/n^{6/5}$.}
 \label{fig:b35124}
\end{minipage}
\end{figure}

 This sequence of lower bounds can also be extrapolated against $1/n^{6/5},$ as shown in Fig.~\ref{fig:b35124} and results in an estimate
of $\mu$ consistent with that obtained by the direct analysis of the original, extended series.

\subsection{ Av(42351)}

This series was known up to, and including, terms of order $x^{15},$ given in the OEIS as sequence A256200. We have extended the series by twelve further terms. We have also extended the sequence of ratios, and the coefficient sequence, by 100 further terms, as described above.
The analysis parallels those cases already discussed.  

\begin{figure}[h!] 
\begin{minipage}[t]{0.45\textwidth} 
\centerline{\includegraphics[width=\textwidth,angle=0]{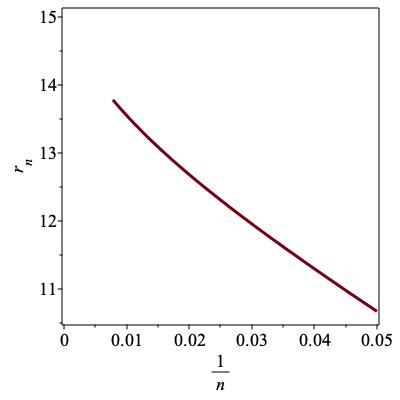} }
\caption{Ratios of $Av(42351)$ vs. $1/n.$} \label{fig:r142351}
\end{minipage}
\hspace{0.05\textwidth}
\begin{minipage}[t]{0.45\textwidth} 
\centerline{\includegraphics[width=\textwidth,angle=0]{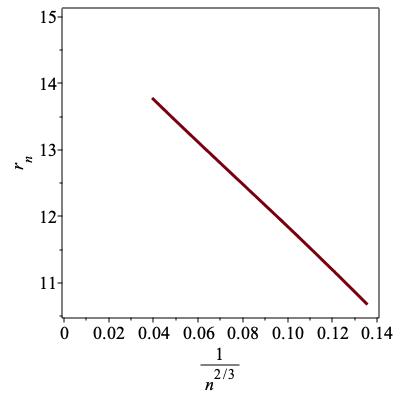}}
\caption{Ratios of $Av(42351)$ vs. $1/n^{2/3}.$} \label{fig:r242351}
\end{minipage}
\end{figure}

The ratio plots of the known series displays curvature when plotted against $1/n,$ see Fig.~\ref{fig:r142351}, and this curvature persists when we include the extrapolated ratios. From the permutation pattern 41325 (which is in the same Wilf class as this permutation), we see the characteristic 1324 pattern that is a hallmark of a stretched exponential term.

The  curvature disappears when the ratios are plotted against $1/n^{2/3},$ see Fig.~\ref{fig:r242351}. This implies a singularity structure with a stretched exponential term. As the pattern $35241$ contains a sub-pattern qualitatively similar to $4231,$ believed to have stretched-exponential asymptotics, this is unsurprising, and, indeed, expected.

The ratio plot against $1/n^{2/3}$ is seen to be visually linear, and extrapolates to $\mu=15.1 \pm 0.1.$ Note that this linearity implies a stretched-exponential term with exponent $\sigma=1/3.$ Linear intercepts, which eliminate the competing $O(1/n)$ term, allow this estimate to be sharpened to $\mu=15.10 \pm 0.05.$

\begin{figure}[h!] 
\begin{minipage}[t]{0.45\textwidth} 
\centerline{\includegraphics[width=\textwidth,angle=0]{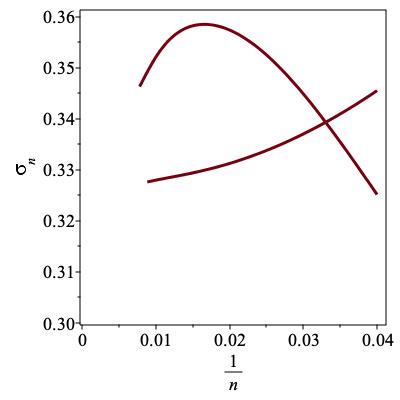} }
\caption{Estimators of $\sigma$ for $Av(42351)$ vs. $1/n$ assuming $\mu=15.0$} \label{fig:s142351}
\end{minipage}
\hspace{0.05\textwidth}
\begin{minipage}[t]{0.45\textwidth} 
\centerline{\includegraphics[width=\textwidth,angle=0]{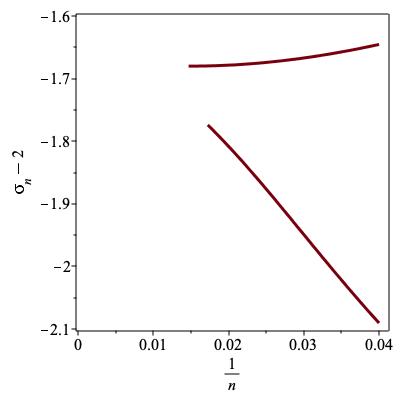}}
\caption{Estimators of $2-\sigma$ for $Av(42351)$ without assuming $\mu.$} \label{fig:s242351}
\end{minipage}
\end{figure}

Assuming $\mu = 15.10,$ as shown above we can estimate the value of $\sigma$ using the two methods we have previously discussed. The estimates of $\sigma$ from these two methods are shown Fig.~\ref{fig:s142351}, and are consistent with our conjectured value, $\sigma$ is $1/3.$ We can also estimate $\sigma$ less precisely without assuming the value of $\mu.$ Again we utilise the two distinct methods discussed above for this estimate. These estimates of $\sigma$ are shown Fig.~\ref{fig:s242351}, and are also consistent with our conjectured value, $\sigma=1/3.$ 

From the expected asymptotic behaviour of the ratios in this case, shown in eqn. (\ref{eq:third}), we can try and estimate the subdominant growth constant $\mu_1$ and exponent $g.$ 
We can use eqn.(\ref{eqn:mu1}) and extrapolate against $1/n,$ and estimate the quantity $\sigma \log{\mu_1}$ as $-6.75 \pm 0.45,$ where we assume $\sigma=1/3,$ and the uncertainty quoted arises from the uncertainty associate with the estimate of $\mu.$ This implies $\mu_1 = 0.0012 \pm 0.0006.$ We emphasis that this is a rather imprecise estimate, and should be considered to be more an order-of-magnitude estimate than anything more. 

We obtained an alternative estimate of $\mu_1$ by direct fitting to the unknown parameters $C,$ $\log{\mu_1}$ and $g$ in
eqn (\ref{logcan1}), assuming the central estimate of $\mu.$ In this way we estimate $g = -1.25 \pm 0.25,$ and obtained a similar estimate of $\mu_1$ to that just quoted. The estimate of $C$ is very sensitive to the estimate of $\mu,$ and we can only conclude $1 < C < 2.$

We therefore conclude that $s_n(42351)\sim C \cdot \mu^n \cdot \mu_1^{n^{\sigma}}  \cdot n^g,$ with $\mu = 15.10 \pm 0.05, \,\, \sigma = 1/3,\,\, \mu_1 =0.0012 \pm 0.0006, \,\, g = -1.25 \pm 0.25,$  and $1 < C < 2.$ 

\subsubsection{Lower bounds}
Assuming that the coefficients form a Stieltjes moment sequence, the consequent log-convexity gives the bound $\mu(42351) \ge 11.4884.$ If the 100 predicted coefficients are believed, this improves the bound to $\mu(42351) \ge 13.7711.$

Constructing the continued fraction representation from the exact coefficients, (see Theorem 1 above), as  above, then one finds $\mu(42351) \ge 14.0314.$

\begin{figure}[h!] 
\centerline{\includegraphics[width=3in,angle=0]{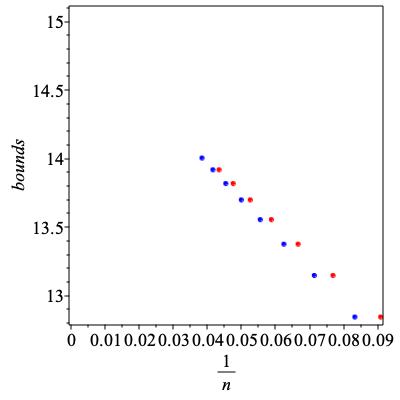} }
\caption{Stieltjes bounds for $Av(42351)$ vs. $1/n^{6/5}.$ } \label{fig:b42351}
\end{figure}

 This sequence of lower bounds can also be extrapolated against $1/n,$ see Fig.~\ref{fig:b42351} and results in an estimate
of $\mu$ consistent with that obtained by the direct analysis of the original, extended series.

\subsection{ Av(42315)}

This series was  known up to, and including, terms of order $x^{14},$ given in the OEIS as sequence A256206. We have extended the series by twelve further terms. We have also extended the sequence of ratios (and coefficients) by 200 further ratios (and coefficients), as described above. 

\begin{figure}[h!] 
\begin{minipage}[t]{0.45\textwidth} 
\centerline{\includegraphics[width=\textwidth,angle=0]{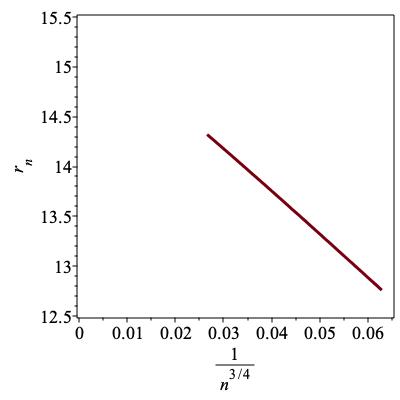} }
\caption{Ratios of $Av(42315)$ vs. $1/n^{3/4}$.} \label{fig:r42315}
\end{minipage}
\hspace{0.05\textwidth}
\begin{minipage}[t]{0.45\textwidth} 
\centerline{\includegraphics[width=\textwidth,angle=0]{A42351bounds.jpg}}
\caption{Stieltjes bounds for $Av(42315)$ vs. $1/n^{6/5}$.} \label{fig:b42315}
\end{minipage}
\end{figure}

The ratios, when plotted against $1/n,$ display curvature, which disappears when the ratios are plotted against $1/n^{3/4},$ see  Fig.~\ref{fig:r42315}.  It appears to be linear, and extrapolates to $\mu=15.4 \pm 0.2.$ This linearity implies a stretched-exponential term with exponent $\sigma=1/4.$ Linear intercepts, which eliminate the competing $O(1/n)$ term, allow this estimate to be sharpened to $\mu=15.40 \pm 0.1.$

Assuming $\sigma=1/4,$ one can fit sequences of four successive ratios to the asymptotic form given by eqn. (\ref{eqn:ratfit}), and also fit four successive coefficients to the coefficient asymptotic form eqn. (\ref{logcan1}). 

From the fit to the ratios, we find $\mu=15.40 \pm 0.04,$ and $\sigma \log{\mu_1} = -2.3 \pm 0.1,$ so $\mu_1 = 0.0001 \pm 0.00005.$ 

From the fit to the logarithm of the coefficients we find $\log{C} = 4.5 \pm 0.2,$ so $C = 90 \pm 16,$ and $\log{\mu} \approx 2.736,$ so $\mu \approx 15.42.$

We combine these various estimates of $\mu$ to give  the final estimate  $\mu = 15.40 \pm 0.04.$ 

\begin{figure}[h!] 
\begin{minipage}[t]{0.45\textwidth} 
\centerline{\includegraphics[width=\textwidth,angle=0]{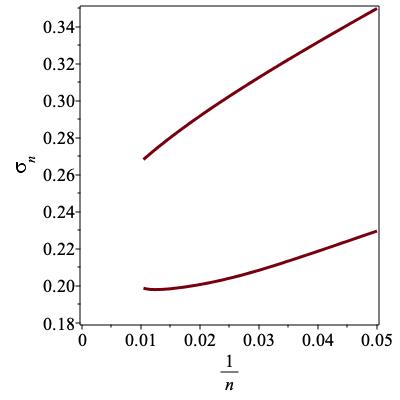}}
\caption{Estimate of exponent $\sigma$ vs. $1/n,$ assuming $\mu.$ First method, upper curve, second method, lower curve.}
\label{fig:sig42315}
\end{minipage}
\hspace{0.05\textwidth}
\begin{minipage}[t]{0.45\textwidth} 
\centerline{\includegraphics[width=\textwidth,angle=0]{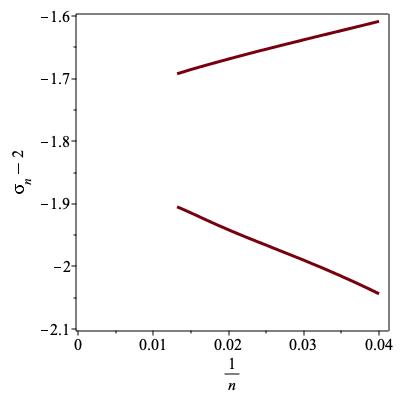} }
\caption{Estimate of exponent $\sigma-2$ vs. $1/n$ without knowing $\mu.$ First method, upper curve, second method, lower curve.}
\label{fig:sig242315}
\end{minipage}
\end{figure}

We directly estimate $\sigma$ assuming the value of $\mu$ by two methods. The first estimator is given by eqn (\ref{eqn:sig1}), the second by calculating the gradient of the log-log plot obtained from eqn(\ref{eq:sigma2}). These estimators are shown in Fig.~\ref{fig:sig42315}.
One can also estimate the value of $\sigma$ without any estimate of $\mu,$ by calculating the local gradient of the log-log plots of the two estimators given by eqn. (\ref{eq:sig1}) and eqn. (\ref{eq:sig2}). These estimators are shown in Fig.~\ref{fig:sig242315}. Both pairs of plots give results consistent with our estimate $\sigma = 1/4.$

We therefore conclude that $s_n(42315)\sim C \cdot \mu^n \cdot \mu_1^{n^{1/4}}  \cdot n^g,$ with $\mu = 15.40 \pm 0.04, \,\, \sigma = 1/4, \,\, \mu_1 = 0.0001 \pm 0.00005,$ and $C = 90 \pm 16.$ We give no estimate of  the exponent $g.$ 

\subsubsection{Lower bounds}
Assuming that the coefficients form a Stieltjes moment sequence, the consequent log-convexity gives the bound $\mu(42315) \ge 11.8117.$ If the 200 predicted coefficients are believed, this improves the bound to $\mu(42315) \ge 14.7160.$

Constructing the continued fraction representation from the exact coefficients, (see Theorem 1 above), as in the previous PAP, one finds  $\mu(42315) \ge 14.5633.$

 This sequence of lower bounds can  be extrapolated against $1/n^{6/5},$ shown in Fig.~\ref{fig:b42315} and results in an estimate
of $\mu \approx 15.4$ consistent with that obtained by the direct analysis of the original, extended series.

\subsection{ Av(35241)}

This series was  known up to, and including, terms of order $x^{14},$ given in the OEIS as sequence A256201. We have extended the series by twelve further terms. We have also extended the sequence of ratios approximately by 100 further coefficients and ratios, as described above. 

\begin{figure}[h!] 
\begin{minipage}[t]{0.45\textwidth} 
\centerline{\includegraphics[width=\textwidth,angle=0]{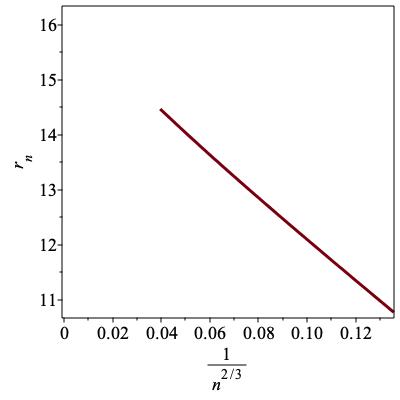} }
\caption{Ratios of $Av(35241)$ vs. $1/n^{2/3}$.} \label{fig:r35241}
\end{minipage}
\hspace{0.05\textwidth}
\begin{minipage}[t]{0.45\textwidth} 
\centerline{\includegraphics[width=\textwidth,angle=0]{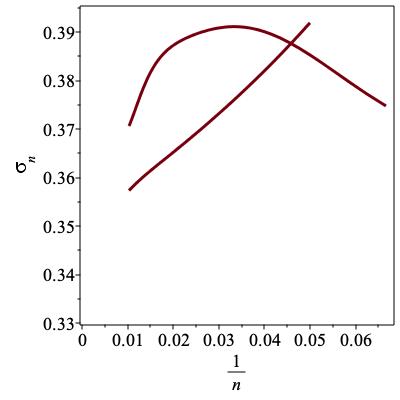}}
\caption{Estimators of $1-\sigma$ for $Av(35241)$ vs. $1/n$.} \label{fig:s35241}
\end{minipage}
\end{figure}

The ratios, when plotted against $1/n,$ display curvature, which disappears when the ratios are plotted against $1/n^{2/3},$ see Fig.~\ref{fig:r35241}. This implies a singularity structure with a stretched exponential term. As the pattern $35241$ contains a sub-pattern qualitatively similar to $4231,$ believed to have stretched-exponential asymptotics, this is unsurprising, and, indeed, expected.

The ratio plot against $1/n^{2/3}$ is seen to be visually linear, and extrapolates to $\mu=16.2 \pm 0.2.$ This linearity implies a stretched-exponential term with exponent $\sigma=1/3.$ Linear intercepts, which eliminate the competing $O(1/n)$ term, allow this estimate to be sharpened to $\mu=16.2 \pm 0.1.$

Assuming $\sigma=1/3,$ one can fit sequences of four successive ratios to the asymptotic form given by eqn. (\ref{eqn:ratfit}), and also fit four successive coefficients to the asymptotic form eqn. (\ref{logcan1}). 

From the fit to the ratios, we find $\mu=16.20 \pm 0.05,$ and $\sigma \log{\mu_1} = -3.1 \pm 0.2,$ so $\mu_1 = 0.00009 \pm 0.00004.$ 

From the fit to the logarithm of the coefficients we find $\log{C} = 2.9 \pm 0.4,$ so $C = 18 \pm 9.$ 

We combine these various estimates of $\mu$ to give  the final estimate  $\mu = 16.20 \pm 0.05.$ 

We can directly estimate $\sigma$ assuming the value of $\mu$ by two methods, as described above. The first estimator is given by eqn. (\ref{eqn:sig1}), the second by calculating the gradient of the log-log plot obtained from eqn(\ref{eq:sigma2}). These estimators are shown in Fig.~\ref{fig:s35241}.
One can also estimate the value of $\sigma$ without any estimate of $\mu,$ by calculating the local gradient of the log-log plots of the two estimators given by eqn. (\ref{eq:sig1}) and eqn. (\ref{eq:sig2}). These estimators are shown in Fig.~\ref{fig:sig35241}. Both pairs of plots give results consistent with our estimate $\sigma = 1/3.$

We therefore conclude that $s_n(35241)\sim C \cdot \mu^n \cdot \mu_1^{n^{\sigma}}  \cdot n^g,$ with $\mu = 16.20 \pm 0.05, \,\, \sigma = 1/3,\,\, \mu_1 \approx 0.00009,$ and $C = 18 \pm 9.$  We give no estimate of the  exponent $g.$

\subsubsection{Lower bounds}
Assuming that the coefficients form a Stieltjes moment sequence, the consequent log-convexity gives the bound $\mu(35241) \ge 11.6779.$ If the 100 predicted coefficients are believed, this improves the bound to $\mu(35241) \ge 14.4634.$

Constructing the continued fraction representation from the exact coefficients, (see Theorem 1 above), as above, then one finds  $\mu(35241) \ge 14.6253.$

\begin{figure}[h!] 
\begin{minipage}[t]{0.45\textwidth} 
\centerline{\includegraphics[width=\textwidth,angle=0]{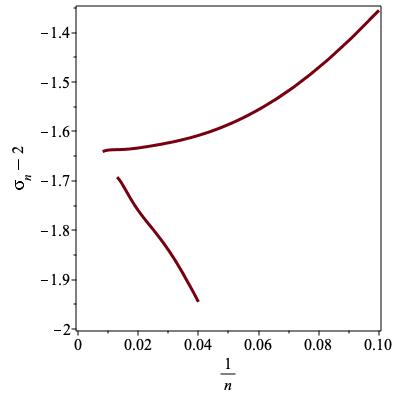}}
\caption{Estimate of exponent $\sigma$ vs. $1/n,$ without $\mu.$ First method, upper curve, second method, lower curve.}
\label{fig:sig35241}
\end{minipage}
\hspace{0.05\textwidth}
\begin{minipage}[t]{0.45\textwidth} 
\centerline{\includegraphics[width=\textwidth,angle=0]{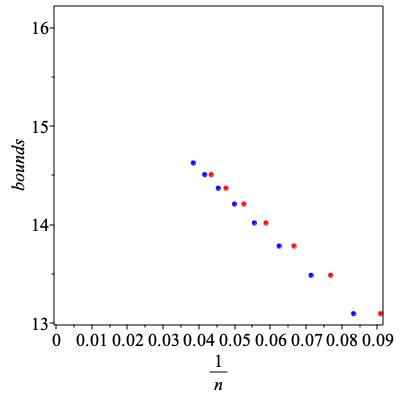} }
\caption{Stieltjes bounds on $\mu$ for $Av(35241).$}
\label{fig:b35241}
\end{minipage}
\end{figure}

 This sequence of lower bounds can  be extrapolated against $1/n,$ shown in Fig.~\ref{fig:b35241}, and results in an estimate
of $\mu \approx 16,$ consistent with that obtained by the direct analysis of the original, extended series.

\subsection{ Av(53241)}
This series was  known up to, and including, terms of order $x^{14},$ given in the OEIS as sequence A256202. We have extended the series by twelve further terms. We have also extended the sequence of ratios, and coefficients, by 100  approximate further terms as described above. 

\begin{figure}[h!] 
\begin{minipage}[t]{0.45\textwidth} 
\centerline{\includegraphics[width=\textwidth,angle=0]{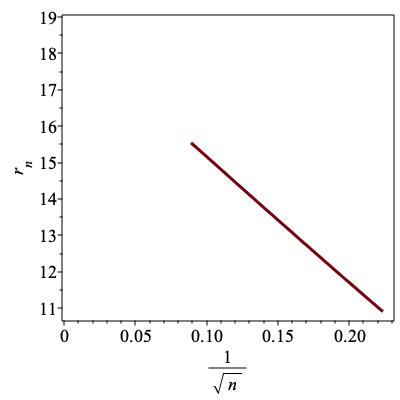} }
\caption{Ratios of $Av(53241)$ vs. $1/\sqrt{n}$.} \label{fig:r53241a}
\end{minipage}
\hspace{0.05\textwidth}
\begin{minipage}[t]{0.45\textwidth} 
\centerline{\includegraphics[width=\textwidth,angle=0]{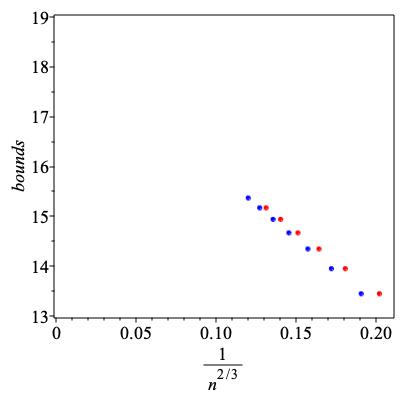} }
\caption{Stieltjes bounds on $\mu$ for $Av(53241).$}
\label{fig:b53241}
\end{minipage}
\end{figure}

The ratios, when plotted against $1/n,$ display curvature which disappears when the ratios are plotted against $1/n^{1/2}.$ This implies a singularity structure with a stretched exponential term. The corresponding ratio plot is shown in Fig.~\ref{fig:r53241a}. It is seen to be visually linear, and extrapolates to $\mu=18.7 \pm 0.2.$ Note that this linearity implies a stretched-exponential term with exponent $1/2.$ Linear intercepts, which eliminate the competing $O(1/n)$ term, allow this estimate to be sharpened to $\mu=18.65 \pm 0.1.$
A variety of other methods were also employed to estimate $\mu.$ These were all consistent, and gave slightly more precision. We combine these various estimates of $\mu$ to give  the final estimate  $\mu = 18.66 \pm 0.05.$ 

Assuming $\sigma=1/2,$ one can fit sequences of four successive ratios to the asymptotic form given by eqn. (\ref{eqn:ratfit}), and also fit four successive coefficients to the asymptotic form eqn. (\ref{logcan1}). 

From the fit to the ratios, we find $\mu=18.66 \pm 0.05,$ and $\sigma \log{\mu_1} = -1.8 \pm 0.1,$ so $\mu_1 = 0.027 \pm 0.006.$ 

From the fit to the logarithm of the coefficients we were unable to usefully estimate $\log{C}.$ Both these methods gave an estimated result for the exponent $g$ in the vicinity of $-4.$ We expect $g=-4 \pm 1.$

We can directly estimate $\sigma$ assuming the value of $\mu$ by two methods, as described above. The first estimator is given by eqn. (\ref{eqn:sig1}), the second by calculating the gradient of the log-log plot obtained from eqn. (\ref{eq:sigma2}). 
One can also estimate the value of $\sigma$ without any estimate of $\mu,$ by calculating the local gradient of the log-log plots of the two estimators given by eqn. (\ref{eq:sig1}) and eqn. (\ref{eq:sig2}). 
Both pairs of plots give results totally consistent with our previous estimate $\sigma = 1/2.$

We therefore conclude that $s_n(53241)\sim C \cdot \mu^n \cdot \mu_1^{n^{\sigma}}  \cdot n^g,$ with $\mu = 18.66 \pm 0.05, \,\, \sigma = 1/2,\,\, \mu_1 = 0.027 \pm 0.006,$ and $g = -4 \pm 1.$ We give no estimate of the  amplitude $C.$

\subsubsection{Lower bounds}
Assuming that the coefficients form a Stieltjes moment sequence, the consequent log-convexity gives the bound $\mu(53241) \ge 11.8339.$ If the 100 predicted coefficients are believed, this improves the bound to $\mu(53241) \ge 15.5411.$

Constructing the continued fraction representation from the exact coefficients, (see Theorem 1 above), as in the previous PAP, then  one finds  $\mu(53241) \ge 15.4445.$

 This sequence of lower bounds can also be extrapolated against $1/n^{2/3},$ shown in Fig.~\ref{fig:b53241}, and results in an estimate
of  $\mu \approx 18.7,$ consistent with that obtained by the direct analysis of the original, extended series.

\subsection{ Av(53421)}
This series was  known up to, and including, terms of order $x^{14},$ given in the OEIS as sequence A256207. We have extended the series by twelve further terms. We have also extended the sequence of ratios, and coefficients, by 200 further terms, as described above. 

\begin{figure}[h!] 
\begin{minipage}[t]{0.45\textwidth} 
\centerline{\includegraphics[width=\textwidth,angle=0]{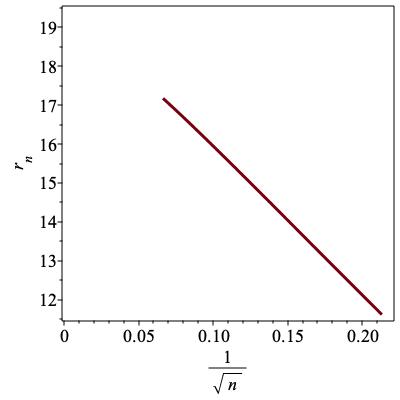} }
\caption{Ratios of $Av(53421)$ vs. $1/\sqrt{n}$.}
\label{fig:r53421}
\end{minipage}
\hspace{0.05\textwidth}
\begin{minipage}[t]{0.45\textwidth} 
\centerline{\includegraphics[width=\textwidth,angle=0]{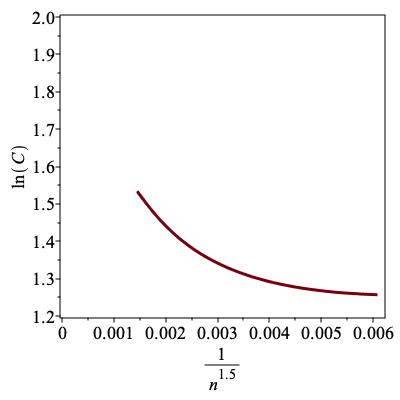} }
\caption{Estimators of $\log(C)$ vs. $1/n^{5/2}$.}
\label{fig:B53421}
\end{minipage}
\end{figure}

The ratios, when plotted against $1/n,$ display curvature, which disappears when the ratios are plotted against $1/n^{1/2}.$ This implies a singularity structure with a stretched exponential term. The ratio plot against $1/\sqrt{n}$ is shown in Fig.~\ref{fig:r53421}. It appears to be linear, and extrapolates to $\mu=19.4 \pm 0.2.$  This linearity implies a stretched-exponential term with exponent $\sigma=1/2.$ 

However, we can get a much more precise estimate of $\mu$ in this case, as
B\'ona in \cite{B07} points out that, using techniques from \cite{B05} one can prove that 
$\mu(53421)=(1+\sqrt{\mu(Av(1324)})^2 \approx 19.4092.$  See also \cite{APV19}. Here we have used the estimate $\mu(Av(1324)) \approx 11.598 $ given in \cite{CGZ18}. In that paper, it is pointed out that $9+3\sqrt{3}/2 =11.598\ldots.$ If that surd is the exact value of $\mu(Av(1324)),$ it would follow that $\mu(53421)=10+\sqrt{36 +6\sqrt{3}} + 3\sqrt{3}/2 = 19.409265890507264\ldots.$ For our subsequent analysis, it is immaterial which of these two estimates of $\mu$ we use.

Given this more precise estimate of $\mu$ than we have had in most other cases, we can get useful estimates of the other critical parameters.
We fit to $$s_n(53421)=c_n \sim C\cdot \mu^n \cdot \mu_1^{\sqrt{n}} \cdot n^g.$$
Then
$$\log c_n-n\log{\mu}  \sim \log C + \sqrt{n} \log \mu_1 + g \log{n}.$$
We can fit successive triples of coefficients to this equation, giving a linear system the solution of which gives estimators of the three unknowns, $C,$ $\log{\mu_1}$ and $g.$
The results of these fits are shown in Figs. \ref{fig:B53421}, \ref{fig:mu153421}, \ref{fig:g53421} below.
Visually extrapolating these gives the estimates $\log{C} = 2 \pm 0.3,$ $\log {\mu_1} = -3.13 \pm 0.03,$ and $g = -4.0 \pm 0.5.$ This gives
$C = 7.4 \pm 2,$ $\mu_1 = 0.044 \pm 0.002,$ and $g = -4.0 \pm 0.5.$

 \begin{figure}[h!] 
\begin{minipage}[t]{0.45\textwidth} 
\centerline{\includegraphics[width=\textwidth,angle=0]{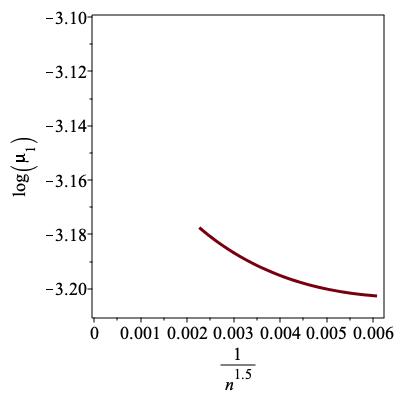}}
\caption{Estimators of $\log(\mu_1)$ vs. $1/n^{3/2}$.}
\label{fig:mu153421}
\end{minipage}
\hspace{0.05\textwidth}
\begin{minipage}[t]{0.45\textwidth} 
\centerline{\includegraphics[width=\textwidth,angle=0]{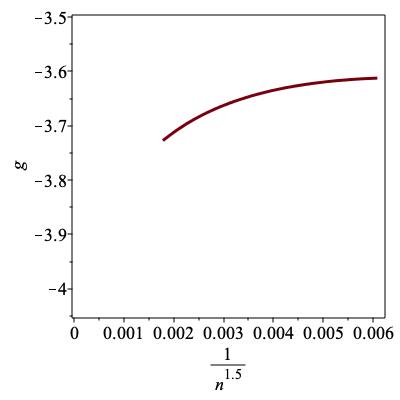} }
\caption{Estimators of $g$ vs. $1/n^{3/2}$.}
\label{fig:g53421}
\end{minipage}
\end{figure}

We therefore conclude that $s_n(53421)\sim C \cdot \mu^n \cdot \mu_1^{n^{\sigma}}  \cdot n^g,$ with $\mu \approx 19.4092, \,\, \sigma = 1/2,$ $C = 7.4 \pm 2,$ $\mu_1 = 0.044 \pm 0.002,$ and $g = -4.0 \pm 0.5.$

 \subsubsection{Lower bounds}
Assuming that the coefficients form a Stieltjes moment sequence, the consequent log-convexity gives the bound $\mu(53421) \ge 12.4079.$ If the 200 predicted coefficients are believed, this improves the bound to $\mu(53421) \ge 17.1769.$

Constructing the continued fraction representation from the exact coefficients, (see Theorem 1 above), as in the previous PAP, then  one finds  $\mu(53421) \ge 16.3053.$

 This sequence of lower bounds can also be extrapolated against $1/n^{2/3},$ and results in an estimate
of  $\mu \approx 18.9,$ a little lower than the known result.

\subsection{ Av(52341)}
This series was  known up to, and including, terms of order $x^{14},$ given in the OEIS as sequence A256208. We have extended the series by nine further terms. This was, computationally, the most demanding series to generate. With 2TB of memory, we could only get to O$(x^{23}),$ with 15 hours of computing time, whereas with most other patterns we could get to O$(x^{26})$ with only 1TB memory and a similar amount of computing time. We have also extended the sequence of ratios and coefficients by 50 further terms as described above. 

\begin{figure}[h!] 
\begin{minipage}[t]{0.45\textwidth} 
\centerline{\includegraphics[width=\textwidth,angle=0]{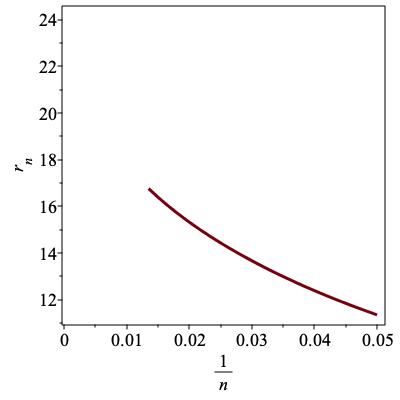} }
\caption{Ratios of $Av(53421)$ vs. $1/n$.}
\label{fig:r52341}
\end{minipage}
\hspace{0.05\textwidth}
\begin{minipage}[t]{0.45\textwidth} 
\centerline{\includegraphics[width=\textwidth,angle=0]{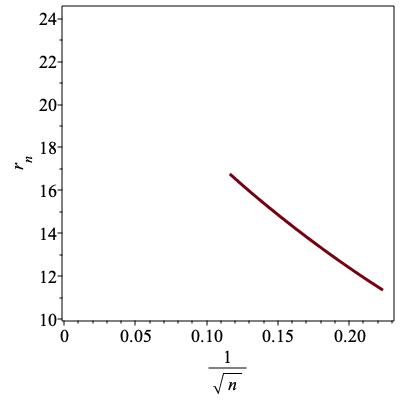} }
\caption{Ratios of $Av(53421)$ vs. $1/\sqrt{n}$.}
\label{fig:r252341}
\end{minipage}
\end{figure}

\begin{figure}[h!] 
\begin{minipage}[t]{0.45\textwidth} 
\centerline{\includegraphics[width=\textwidth,angle=0]{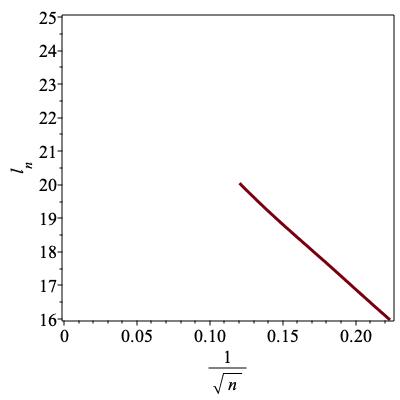} }
\caption{Linear intercepts of $Av(53421)$ vs. $1/\sqrt{n}$.}
\label{fig:l52341}
\end{minipage}
\hspace{0.05\textwidth}
\begin{minipage}[t]{0.45\textwidth} 
\centerline{\includegraphics[width=\textwidth,angle=0]{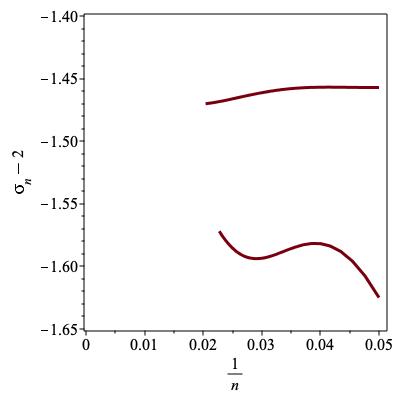} }
\caption{Estimators of $\sigma$ for $Av(53421)$ vs. $1/n$.}
\label{fig:sig52341}
\end{minipage}
\end{figure}

The ratios, when plotted against $1/n,$ see Fig.~\ref{fig:r52341}, display considerable curvature, which is substantially reduced when the ratios are plotted against $1/n^{1/2},$ shown in Fig.~\ref{fig:r252341}. This again suggests a singularity structure with a stretched exponential term, similar to that observed for $Av(4231)$ PAPs \cite{CGZ18}. This is not surprising, as this pattern contains the pattern $4231.$ Extrapolating the ratios, we estimate $\mu \approx 24.6 \pm 0.5.$ Next, we calculated the linear intercepts $l_n \equiv n\cdot r_n - (n-1) \cdot r_{n-1}$ which eliminate the O$(1/n)$ term in the ratios. The results are shown in Fig.~\ref{fig:l52341}, and which we extrapolate to $\mu \approx 24.8 \pm 0.5.$ 

Using the methods discussed in the analysis of previous patterns, we have estimated the value of the stretched-exponential exponent $\sigma.$ We show in Fig.~\ref{fig:sig52341} estimators of $\sigma$ obtained without assuming the value of $\mu.$ It can be seen that they are quite consistent with the value $1/2,$ so we feel confident suggesting that $\sigma = 1/2$ exactly.

\begin{figure}[h!] 
\begin{minipage}[t]{0.45\textwidth} 
\centerline{\includegraphics[width=\textwidth,angle=0]{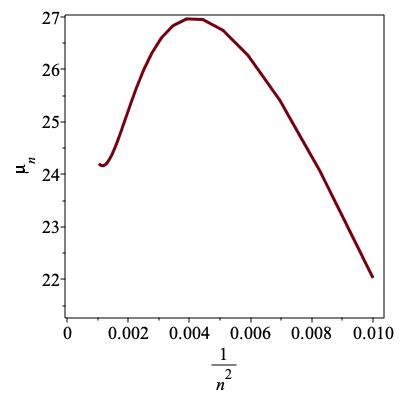}}
\caption{Estimators of $\mu$  for $Av(52341)$ vs. $1/n^{2}$.}
\label{fig:mu52341}
\end{minipage}
\hspace{0.05\textwidth}
\begin{minipage}[t]{0.45\textwidth} 
\centerline{\includegraphics[width=\textwidth,angle=0]{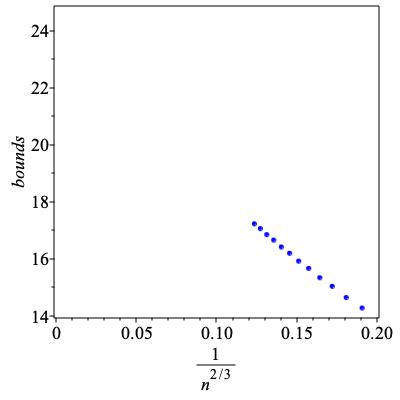} }
\caption{Bounds on $\mu$ for $Av(52341)$ vs. $1/n^2/3$.}
\label{fig:b52341}
\end{minipage}
\end{figure}

Assuming $\sigma=1/2$ have also estimated the value of $\mu$ by fitting successive quadruples of ratios, $r_{k-2},\, r_{k-1},\, r_k,\, r_{k+1}$ to the expected asymptotic form (\ref{eq:half}), which gives estimators of $\mu$ as shown in Fig.~\ref{fig:mu52341}, which while difficult to extrapolate, is not inconsistent with previous estimates. This fit also gives estimators of the growth constant $\mu_1,$ more precisely of $\mu \cdot \log{\mu_1}/2,$ which extrapolates to a value around $-80,$ from which one concludes that $\mu_1 \approx 0.0016.$ We can also estimate the growth constant $\mu_1$ by plotting $(r_n/\mu-1) \sim \log(\mu_1)/2$ against $1/\sqrt{n},$ which extrapolates to $-3.3 \pm 0.1,$ so that $\mu_1 = 0.0014 \pm 0.0003.$

We therefore conclude that $s_n(52341)\sim C \cdot \mu^n \cdot \mu_1^{\sqrt{n}}  \cdot n^g,$ with $\mu = 24.8 \pm 0.5, \,\, \sigma = 1/2, \,\ \mu_1 \approx 0.0015. $
We give no estimate of  the exponent $g$ or the amplitude $C.$

\subsubsection{Lower bounds}
Assuming that the coefficients form a Stieltjes moment sequence, the consequent log-convexity gives the bound $\mu(52341) \ge 12.1999.$ If the 50 predicted coefficients are believed, this improves the bound to $\mu(52341) \ge 16.7641.$

If one constructs the continued fraction representation from the exact coefficients, (see Theorem 1 above), as in the previous PAP, then by the same construction one finds  $\mu(52341) \ge 15.7310.$

 This sequence of lower bounds can also be extrapolated against $1/n^{2/3},$ see Fig.~\ref{fig:b52341} and results in an estimate
of $\mu$ consistent with that obtained by the direct analysis of the original, extended series.


\end{document}